\documentclass{amsart}

\usepackage{amsmath,amscd,amssymb,amsthm}
\usepackage{enumerate,mathrsfs}
\usepackage{geometry}
\usepackage{hyperref}
\usepackage{url}
\usepackage{ifthen}
\usepackage[all]{xy} \xyoption{2cell} \UseAllTwocells
\usepackage[mathscr]{euscript}
\usepackage[T1]{fontenc}
\usepackage[utf8]{inputenc}   
\usepackage{refcount}
\usepackage{zref-user}
\usepackage{zref-abspage}

\numberwithin{equation}{section}

\theoremstyle{plain}
\newtheorem{thm_}[equation]{Theorem}
\newtheorem{lemma_}[equation]{Lemma}
\newtheorem{prop_}[equation]{Proposition}
\newtheorem{cor_}[equation]{Corollary}
\newtheorem{eg_}[equation]{Example}

\theoremstyle{definition}
\newtheorem{thmu_}[equation]{Theorem}
\newtheorem*{thmus_}{Theorem}
\newtheorem{propu_}[equation]{Proposition}
\newtheorem*{propus_}{Proposition}
\newtheorem{coru_}[equation]{Corollary}
\newtheorem*{corus_}{Corollary}
\newtheorem{lemu_}[equation]{Lemma}
\newtheorem*{lemus_}{Lemma}
\newtheorem{egu_}[equation]{Example}
\newtheorem*{egus_}{Example}
\newtheorem{def_}[equation]{Definition}
\newtheorem*{defs_}{Definition}
\newtheorem{rk_}[equation]{Remark}
\newtheorem*{rks_}{Remark}
\newtheorem{ex_}[equation]{Remark}
\newtheorem*{exs_}{Remark}
\newtheorem{constr_}[equation]{Construction}
\newtheorem*{constrs_}{Construction}
\newtheorem{nota_}[equation]{Notation}
\newtheorem*{notas_}{Notation}

\newcommand{\thm}[1]{\begin{thm_}#1\end{thm_}}
\newcommand{\thmu}[1]{\begin{thmu_}#1\end{thmu_}}

\newcommand{\lemm}[1]{\begin{lemma_}#1\end{lemma_}}

\newcommand{\egu}[1]{\begin{egu_}#1\end{egu_}}

\newcommand{\prop}[1]{\begin{prop_}#1\end{prop_}}

\newcommand{\defi}[1]{\begin{def_}#1\end{def_}}

\newcommand{\rk}[1]{\begin{rk_}#1\end{rk_}}

\newcommand{\cor}[1]{\begin{cor_}#1\end{cor_}}

\newcommand{\pf}[1]{\begin{proof}#1\end{proof}}

\DeclareMathOperator{\SL}{SL}

\DeclareMathOperator{\Gal}{Gal}
\DeclareMathOperator{\Hom}{Hom}

\DeclareMathOperator{\Ext}{Ext}

\DeclareMathOperator{\Spec}{Spec}

\DeclareMathOperator{\Tr}{Tr}
\DeclareMathOperator{\im}{im}

\DeclareMathOperator{\Br}{Br}

\DeclareMathOperator{\Pic}{Pic}

\DeclareMathOperator{\inv}{inv}

\DeclareMathOperator{\coker}{coker}

\DeclareMathOperator{\KD}{KD}
\DeclareMathOperator{\UU}{U}
\DeclareMathOperator{\id}{id}

\DeclareMathOperator{\colim}{colim}

\DeclareMathOperator{\Map}{Map}

\newcommand{\Chp}{{\sC\tu{hp}}}
\newcommand{\bDel}{{\mathbf{\Delta}}}

\newcommand{\NN}{\mathbb N}

\newcommand{\QQ}{\mathbb Q}

\newcommand{\ZZ}{\mathbb Z}
\newcommand{\bfA}{\mathbf A}%
\newcommand{\bfG}{\mathbf G}%

\newcommand{\sD}{\mathscr D}

\newcommand{\sF}{\mathscr F}
\newcommand{\sG}{\mathscr G}
\newcommand{\sH}{\mathscr H}

\newcommand{\sO}{\mathscr O}

\newcommand{\sS}{\mathscr S}

\newcommand{\sX}{\mathscr X}
\newcommand{\sY}{\mathscr Y}

\newcommand{\cA}{\mathcal A}%
\newcommand{\cB}{\mathcal B}%
\newcommand{\cC}{\mathcal C}%
\newcommand{\cD}{\mathcal D}

\newcommand{\cO}{\mathcal O}%

\newcommand{\cT}{\mathcal T}
\newcommand{\cU}{\mathcal U}

\newcommand{\cW}{\mathcal W}
\newcommand{\cX}{\mathcal X}
\newcommand{\cY}{\mathcal Y}
\newcommand{\cZ}{\mathcal Z}%

\newcommand{\s}{\sigma}

\newcommand{\La}{\Lambda}
\newcommand{\tm}{\times}%
\newcommand{\otm}{\otimes}
\newcommand{\ol}{\overline}

\newcommand{\bu}{\bullet}

\newcommand{\ra}{\rightarrow}

\newcommand{\Ra}{\Rightarrow}

\newcommand{\xra}{\xrightarrow}

\newcommand{\hra}{\hookrightarrow}

\newcommand{\rla}{\rightleftarrows}

\newcommand{\mpt}{\mapsto}

\newcommand{\os}[2]{\overset{#1}{#2}}
\newcommand{\us}[2]{\underset{#1}{#2}}
\newcommand{\is}[2]{\xymatrix@-4mm{#1 \ar[r]^-{\sim} & #2 }}
\newcommand{\mis}[2]{\xymatrix@-2mm{#1 \ar[r]^-{\sim} & #2 }}

\newcommand{\dra}[4]{\xymatrix@-4mm{#1 \ar@<.5ex>[r]^-{#3} \ar@<-.5ex>[r]_-{#4}& #2 }}  
\newcommand{\era}[5]{\xymatrix@-4mm{#1 \ar[r] &#2 \ar@<.5ex>[r]^-{#4} \ar@<-.5ex>[r]_-{#5}& #3 }}  

\newcommand{\tu}[1]{\text{\upshape #1}}

\newcommand{\SSet}{\tu{Set}}

\newcommand{\AAb}{\tu{Ab}}
\newcommand{\MMod}{\tu{Mod}}
\newcommand{\PPSh}{\tu{PSh}}
\newcommand{\SSh}{\tu{Sh}}
\newcommand{\SShAb}{\tu{ShAb}}

\newcommand{\SShGrp}{\tu{ShGrp}}

\newcommand{\SSch}{\tu{Scheme}}
\newcommand{\SSpc}{\tu{Space}}
\newcommand{\SStk}{\tu{Stack}}

\newcommand{\op}{{op}}

\newcommand{\et}{\tu{\'et}}
\newcommand{\liset}{\tu{lis-\'et}}
\newcommand{\cart}{\tu{cart}}

\newcommand{\fppf}{\tu{fppf}}

\newcommand{\ctf}{\tu{ctf}}

\DeclareFontFamily{U}{wncy}{}
\DeclareFontShape{U}{wncy}{m}{n}{%
   <5>wncyr5%
   <6>wncyr6%
   <7>wncyr7%
   <8>wncyr8%
   <9>wncyr9%
   <10>wncyr10%
   <11>wncyr10%
   <12>wncyr6%
   <14>wncyr7%
   <17>wncyr8%
   <20>wncyr10%
   <25>wncyr10}{}
\DeclareMathAlphabet{\cyrille}{U}{wncy}{m}{n}
\def\Sha{\cyrille X}
\def\Be{\cyrille B}
\newcommand{\Tors}{\mathsf{Tors}}

\newcommand{\eq}[1]{\begin{equation}#1\end{equation}}
\newcommand{\eqn}[1]{\begin{equation*}#1\end{equation*}}
\newcommand{\ga}[1]{\begin{gather}#1\end{gather}}
\newcommand{\gan}[1]{\begin{gather*}#1\end{gather*}}

\newcommand{\aln}[1]{\begin{align*}#1\end{align*}}

\newcommand{\enmt}[1]{\begin{enumerate}#1\end{enumerate}}

\newcommand{\aci}[1]{\ar@{^(->}[#1]|-{/}}
\newcommand{\coaci}[1]{\ar@{_(->}[#1]|-{/}}
\newcommand{\aoi}[1]{\ar@{^(->}[#1]|-{\circ}}
\newcommand{\coaoi}[1]{\ar@{_(->}[#1]|-{\circ}}
\makeatletter
\def\citet@url@sp{https://stacks.math.columbia.edu/}
\def\citet@bib@sp{stacks-project}
\def\citet@url@kd{https://kerodon.net/}
\def\citet@bib@kd{kerodon}
\newcommand{\citet@tag}[2]{\href{#2tag/#1}{#1}}
\newcommand{\citet@taglist}[2]{%
 \def\@citet@e{}%
 \def\@citet@tag@n{0}
 \@for\@citet@tag:=#1\do{%
  \edef\@citet@tag@n{\the\numexpr\@citet@tag@n + 1}%
 }%
 \def\@citet@tags{%
  \def\@citet@tag@i{0}%
  \@for\@citet@tag:=#1\do{%
   \edef\@citet@tag@i{\the\numexpr\@citet@tag@i + 1}%
   \ifthenelse{\@citet@tag@i > 1}{
    \ifthenelse{\@citet@tag@i = \@citet@tag@n}{
     \citet@seplast%
    }{%
     \citet@sep%
    }%
   }{}%
   \citet@entry{\@citet@tag}{#2}
  }%
 }%
 \ifthenelse{\@citet@tag@n > 1}{
  \def\@citet@Tag{Tags}%
 }{%
  \def\@citet@Tag{Tag}%
 }%
 \@citet@Tag~\@citet@tags
}
\newcommand{\citet@sep}{, }
\newcommand{\citet@seplast}{ and }
\newcommand{\citet@entry}[2]{\citet@tag{#1}{#2}}
\let\@old@cite\cite
\renewcommand{\cite}[2][]{%
 \def\@citet@detail{\citet@taglist{#1}{\@citet@url}}%
 \ifthenelse{\equal{#2}{sp}}{%
  \def\@citet@url{\citet@url@sp}%
  \def\@citet@bib{\citet@bib@sp}%
 }{\ifthenelse{\equal{#2}{kd}}{%
  \def\@citet@url{\citet@url@kd}%
  \def\@citet@bib{\citet@bib@kd}%
 }{
  \def\@citet@detail{#1}%
  \def\@citet@bib{#2}%
 }}%
 \ifthenelse{\equal{#1}{}}{%
  \@old@cite{\@citet@bib}%
 }{%
  \@old@cite[\@citet@detail]{\@citet@bib}%
 }%
}
\makeatother

\newcommand{\etale}{{\'etale}}

\newcommand{\Cech}{{\v{C}ech}}
\newcommand{\DM}{{Deligne--Mumford}}
\newcommand{\Grot}{{Grothendieck}}

\newcommand{\Kunn}{{K\"unneth}}

\newcommand{\BM}{{Brauer--Manin}}
\newcommand{\PT}{{Poitou--Tate}}
\newcommand{\TS}{{Tate--Shafarevich}}

\newcommand{\CT}{{Colliot-Th\'el\`ene}}
\newcommand{\Cart}{{Cartesian}}
\newcommand{\Noet}{{Noetherian}}

\newcommand{\adele}{{ad\`ele}}
\newcommand{\adelic}{{ad\`elic}}

\newcommand{\desc}{\tu{desc}}
\newcommand{\conn}{\tu{conn}}

\newcommand{\sdesc}{{2\tu{-}\desc}}

\newcommand{\hdesc}[2]{{#1\tu{-}\desc{\ifthenelse{\equal{#2}{}}{}{_#2}}}}

\newcommand{\St}[1]{\SSh(#1_\tau)}
\newcommand{\XA}{X(\bfA_k)}

\newcommand{\Xk}{X(k)}

\newcommand{\cXk}{\cX(k)}
\newcommand{\cXAk}{\cX(\bfA_k)}
\newcommand{\cYAk}{\cY(\bfA_k)}
\newcommand{\xp}[3]{\langle#1, #2\rangle_\tu{#3}}
\newcommand{\BMp}[2]{\xp{#1}{#2}{BM}}

\renewcommand{\Chp}{\SStk}

\begin{document}
\title[Brauer-{M}anin obstruction on  stacks]
 {The Brauer-{M}anin obstruction on algebraic stacks}

\author[C. Lv]{Chang Lv}
\address{State Key Laboratory of Cyberspace Security Defense\\
Institute of Information Engineering\\
Chinese Academy of Sciences\\
Beijing 100093, P.R. China}
\email{lvchang@amss.ac.cn}

\author[H. Wu]{Han Wu}
	\address{Hubei Key Laboratory of Applied Mathematics,
	Faculty of Mathematics and Statistics,
	Hubei University,
	No. 368, Friendship Avenue, Wuchang District, Wuhan,
	Hubei, 430062, P.R.China.}
\email{wuhan90@mail.ustc.edu.cn}

\date{\today}
\thanks{C. Lv is partially supported by
 National Natural Science Foundation of China NSFC Grant No. 11701552. H. Wu is
 partially funded by Open Foundation of Hubei Key Laboratory of Applied Mathematics (Hubei University) and by the National Natural Science Foundation of	China (NSFC) Grant No. 12301014, and supported by Hubei Provincial Natural Science Foundation of China Grant No. 2023AFB282.
}

\begin{abstract}
For algebraic stacks over number fields,
 we define their {\BM} sets, {\BM} pairings,
 and extend the
  descent theory of   {\CT} and Sansuc.
By extending Sansuc's exact sequence,
 we show the torsionness
 of Brauer groups of stacks that are locally quotients of varieties
 by linear groups.
With mild assumptions, for stacks that are locally quotients or
 {\DM}, we show that the {\BM} set equals the second descent set
 and is contained in the descent set defined by torsors under
 connected linear groups.
For {\BM} sets of these stacks,
 we show the properties such as descent along a torsor and
 product preservation are still correct.
These results extend classical theories of those on varieties.
\end{abstract}
\maketitle

\setcounter{tocdepth}{1}
\tableofcontents

\section{Introduction} \label{intro}
Throughout this paper, $k$ denotes a number field; in particular,
 $\operatorname{char}(k)=0$.

Suppose that  $X$ is a variety over a number field $k$, and
 $\Br X=H_\et^2(X, \bfG_m)$ is the (cohomological) {Brauer--\Grot} group of $X$
 \cite{grothendieck95brauer}.
We have the {\BM}  pairing (see, for example, \cite[8.2]{poonen17rational})
\aln{
\XA \tm \Br X & \ra \QQ/\ZZ, \\
((x_v)_v, A) & \mpt \sum_{v\in\Omega_k}\inv_v A(x_v),
}
 where $\bfA_k$ (resp. $\Omega_k$) is the {\adele} (resp. the set of
 all places) of $k$,
  $\inv_v$ is the invariant map at $v$, and
 $A(x_v) = x_v^*A$  is the image of $A$ under the
  pullback functor $x_v^*: \Br X \ra \Br k_v$.
Under this pairing, the {\BM} set (obstruction) is defined to be the subset
 $\XA^{\Br}\subseteq\XA$ orthogonal to the whole $\Br X$,
 which contains $\Xk$ and  plays an important
 role in some failures of the local-global principle.

In the development of moduli theory, some geometric objects, other than schemes,
 were used to enlarge the category of schemes.
Artin \cite{Ar70},
 Deligne and  Mumford \cite{DM69}
 defined algebraic stacks, and used them to explain some moduli problems.
The arithmetic of these objects is of its own interest.
Let $K$ be a number field and $\cO_K$ be the ring of integers.
Given an algebraic stack $\cX$ over $\cO_K$,
 we say $\cX$ violates the \textit{local-global principle for integral points}
 if $\cX(\cO_{v})\neq\emptyset$ for all $v\in\Omega_K$,
 whereas $\cX(\cO_K)=\emptyset$.
Counterexamples to the local-global principle for integral points  are given in
 \cite{DG95}, \cite{bp22stackc}, \cite{wl23stackc}.
Bhargava and Poonen \cite{bp22stackc} proved that for any stacky curve over
 $\cO_K$ of genus less than $1/2$,
 it satisfies the local-global principle for integral points.
Furthermore, Christensen \cite{christensen20top} extended the {\adelic}
 topology
 to the set $\cXAk$,
 and proved that for any stacky curve over $\cO_K$ of genus less than $1/2$,
  it satisfies strong approximation.
More recently, Santens \cite{santens2022bm-stack}  showed
 that the {\BM} obstruction  to strong approximation is the only one
 for every stacky curve over
 a global field with finite abelian fundamental group.
Loughran and Santens \cite{ls25malle} established a
 new framework for Brauer groups of stacks, 
  and  defined a new notion of the unramified
 Brauer group for algebraic stacks.

In this paper, we consider the {\BM} obstruction on algebraic stacks over
 number fields.
In Section \ref{frame}, we introduce
 points, cohomological obstructions and torsors on algebraic stacks.
In Section \ref{calc}, by extending Sansuc's exact sequence,
 we calculate  the Brauer groups of  stacks that are quotients of varieties
 by connected groups,
 which contain almost all geometric objects appearing in moduli spaces.
As a consequence, we also show that  the Brauer groups
 of some locally quotient stacks are still torsion,
which is a classical result for smooth varieties.
In Section \ref{descentfund},
 for general algebraic stacks,
 we establish the descent theory of {\CT} and Sansuc
 \cite{cs87descente-ii}, and of Harari and Skorobogatov
 \cite{hs13descent},
 and we also define   the {\BM}
 pairing for stacks.
In Section \ref{comparision},
 we show that with mild assumptions, for stacks that are locally quotients or
 {\DM},  the {\BM} obstruction has some relations
 with  some other cohomological obstructions such as obstructions given by
 torsors under connected groups or abelian gerbes, which was originally
 proved by Harari \cite{harari02groupes} for varieties.
In Section \ref{descentbm}, we generalize some results
 of  Cao \cite{cao18approximation} on  descent for the {\BM}
 sets along a torsor under a connected group.
In Section \ref{prod}, we prove that in the same cases as in Section
 \ref{comparision},
 the {\BM} set of product stacks equals the product of the {\BM} sets of
 stacks, extending results of Skorobogatov and Zarhin \cite{sz14product}
  and the first named author
 \cite{lv20brprodge}.

More specifically, let
\eqn{
\cXAk^\conn=\bigcap_{\text{all connected  linear $k$-groups  $G$}}
 \cXAk^{\check H_\fppf^1(-, G)},
}
and
\eqn{
\cXAk^\sdesc=\bigcap_{\text{all commutative  linear $k$-groups  $G$}}
 \cXAk^{H_\et^2(-, G)},
}
 where
 $\cXAk^{\check H_\fppf^1(-, G)}$ (resp. $\cXAk^{H_\et^2(-, G)}$) is defined
  before \eqref{eq_desc} (resp. before Example \ref{eg_bm}).
Then the main results of this paper are described as follows.
\thmu{ [{Theorems \ref{thm_Br=2desc} and \ref{thm_Br=conn}}]
Let $\cX$ be a smooth quasi-separated
 algebraic $k$-stack
 of finite type that is either {\DM} or
Zariski-locally the quotient of a smooth geometrically integral
 $k$-variety by a linear $k$-group. Then
 $\cXAk^{\Br}=\cXAk^{\sdesc} \subseteq \cXAk^{\conn}$.
}

\thmu{ [{Corollary \ref{cor_prod_quo}}]
If $\cX$ and $\cY$ are stacks that are quotients of smooth geometrically integral
 quasi-affine $k$-varieties by connected linear $k$-groups, then
\eqn{
(\cX\tm_k\cY)(\bfA_k)^{\Br}\xra{\sim} \cXAk^{\Br}\tm \cYAk^{\Br}
}}

These results extend classical theories of those on varieties.
The main technique for extending them is cohomological descent
 (see \ref{cohdes}).

\section{Points, obstructions  and torsors on algebraic stacks} \label{frame}
Let us briefly introduce basic notions and facts used in this
 paper, mainly for cohomological obstructions
 to the local-global principles on algebraic stacks, under the framework
 \cite[Section~2]{lv2desc}  but in a slightly modified way.

Given a commutative ring $A$, if there is no confusion,
 we also write $A$ for $\Spec A$.
Let $S$ be a scheme. We consider algebraic $S$-stacks (resp.
 $S$-spaces) or algebraic stacks (resp. spaces) over $S$
 in the sense of \cite[026O]{sp}.
Let $\SStk/S$ denote the $(2, 1)$-category of algebraic stacks over $S$.
For $\cX\in \SStk/S$, a smooth (resp. {\etale}, resp. flat locally of
 finite presentation) surjective
 $1$-morphism $X\ra \cX$ where $X$ is a scheme will be called
 an atlas (resp. {\etale}-atlas, resp. fppf-atlas) for $\cX$.
When $S=\Spec k$, we denote
 $\ol\cX = \cX\tm_k \ol k$, where $\ol k$ is a fixed algebraic closure.

\subsection{Points and obstructions} \label{pt+ob}
Let $k$ be a number field,  $R$ be a commutative $k$-algebra,
 and $\cX$ be an algebraic stack over $k$.
The set of \emph{$R$-points} of $\cX$, denoted by $\cX(R)$, is the set
 of isomorphism classes of objects in the groupoid $\Hom_{\SStk/k}(\Spec R, \cX)$.
Let $\cXAk$ (resp.  $\cX(k)$) be the \emph{{\adelic} points}
 (resp. \emph{rational points}) of $\cX$.
\rk{
Under some conditions, we can also define {\adelic} points explicitly. See
 \eqref{eq_adelic_explicit}.
}

\defi{ \label{defi_stable}
Let $\cC$ be a $(2, 1)$-category, and $\cD$ be an  ordinary
 category.
Let $F: \cC\ra\cD$ be a functor
 from the underlying ordinary category of $\cC$ to $\cD$.
We say that $F$ is \emph{stable} if
 $F$ is a strict $2$-functor from $\cC$ to $\cD$,
 i.e., for any $1$-morphisms $f$ and $g$ in  $\cC$
 that are $2$-isomorphic, we have $F(f)=F(g)$.
}
\rk{
There is a notion of \emph{stable functor} in the literature
 which means ``having a left adjoint on each slice''.
In this paper, we  keep the
  nomenclature in Definition \ref{defi_stable}.
}

For a stable functor $F: (\SStk/k)^\op\ra \SSet$ and $\cX\in \SStk/k$,
 let $q: \Spec \bfA_k\ra \Spec k$ be the map induced by the
 natural inclusion $k\subset \bfA_k$.
Taking $A\in F(\cX)$,
 as in the classical case for varieties
 (cf. \cite[8.1.1]{poonen17rational}),
 one may define
 the \emph{obstruction given by $A$} to be  the subset
 $\cXAk^A$ of $\cXAk$ whose elements are characterized by
\eq{ \label{eq_XAkA}
\cXAk^A\colon =\{x\in \cXAk\mid A(x)\in \im F(q)\},
}
 (which is well-defined by stability of $F$).
 The  \emph{$F$-set} (or \emph{$F$-obstruction})
 is the subset $\cXAk^F$ whose elements are
 characterized by
\eqn{
\cXAk^F\colon=\bigcap_{A\in F(\cX)} \cXAk^A =
 \{x\in \cXAk\mid \im F(x)\subseteq \im F(q)\}.
}
Then $q$ induces
\eqn{
 \cXk\to \cXAk^F\subseteq  \cXAk^A \subseteq  \cXAk.
}

\rk{ If $\cX$ is represented by a separated scheme, then these definitions
 coincide with the classical ones, and the map in the last row is an injection.
But in general, the map  $\cX(k)\ra \cXAk$ is not necessarily injective.
For example, let $G$ be a linear $k$-group.
Let $*$  be the neutral element in  the pointed set
 $BG(\bfA_k) = H_\fppf^1(\bfA_k, G)$.
Then its preimage $\ker(BG(k)\ra BG(\bfA_k))$
 defined by the following {\Cart} diagram
\eqn{\xymatrix{
&\ker(BG(k)\ra BG(\bfA_k))\ar@{=}[ddl] \ar[r] \ar@{^(->}[d] &\{*\}\ar@{^(->}[d] \\
 & BG(k)\ar[r]\ar@{=}[d] & BG(\bfA_k) \ar@{=}[d]\\
\Sha^1(G/k)\ar[r]& H^1(k, G)\ar[r] &\check H_\fppf^1(\bfA_k, G)
}}
  is   the {\TS}   set $\Sha^1(G/k)$, which
 is not necessarily trivial.
}

Anyway, by abusing notation we may let $\cXk$ be its image in $\cXAk$.
Thus we may write $\cXk\subset \cXAk^F$.

\subsection{Cohomology on algebraic stacks and cohomological obstructions}
  \label{coh+Hiob}
Let $\SSch$ be the category of quasi-separated schemes over $\ZZ$.
By a \emph{$k$-variety} we mean an object of $\SSch$ that is separated of
 finite type over $k$.
For any  $\cT\in\SStk/k$,
 let $\cT_\et$ (resp. $\cT_\fppf$) denote the site with topology inherited from
 the big {\etale} site $(\SSch/k)_\et$
 (resp. the big fppf site $(\SSch/k)_\fppf$).
See \cite[06TN,06TP]{sp}.

Let $\PPSh$ (resp. $\SSh$) denote the category of presheaves (resp. sheaves).
Let  $\tau\in\{\et, \fppf\}$ and
 $f: {\cT'}\ra\cT$ be a $1$-morphism in $\SStk/k$.
Then composing with $f$ gives the functor
 $f^*: \PPSh(\cT)\ra \PPSh(\cT')$, which has a right adjoint $f_*$.
Both $f^*$ and $f_*$ send $\tau$-sheaves to $\tau$-sheaves
 \cite[06TS]{sp}, and thus define
 a morphism of topoi
\eqn{
f=(f_*, f^*): \St{\cT'}\ra \St{\cT},
}
 which is also a pair of adjoint functors
 with $\SSh$ replaced by $\SShAb$ the
 category of sheaves of abelian groups (or $\SShGrp$ the
 category of sheaves of groups), and is compatible with the classical case when
 $\cT$ and $\cT'$ are represented by $S$-schemes, cf. \cite[075J]{sp}.

Following \cite[page~4]{lz24enhanced}, an abelian sheaf $\sF$ on
$\cT_\et$ is \emph{Cartesian} if, for every morphism $h:U\ra V$ of schemes
over $\cT$, the comparison map on the small \'etale sites is an isomorphism
\eqn{h^{-1}\sF_V\xra{\sim}\sF_U.}
We write $D_\cart(\cT_\et)$ for the full subcategory of
$D(\cT_\et)$ whose cohomology sheaves are Cartesian.

For $\tau\in\{\et,\fppf\}$, put
\eqn{D(\cT_\tau)=D(\SShAb(\cT_\tau)).}
These are ordinary unbounded derived categories. The exact inverse image
$f^*=f^{-1}$ and $K$-injective resolutions give the adjunctions
\cite[07A5,079P]{sp}
\eq{ \label{eq_biget}
f^*: D(\cT_\et)\rla D(\cT'_\et): Rf_*
}
and
\eqn{
f^*: D(\cT_\fppf)\rla  D(\cT'_\fppf): Rf_*
}
with
\gan{
R\Gamma(\cT_\et, -): D(\cT_\et) \ra D(\AAb), \text{ and} \\
R\Gamma(\cT_\fppf, -): D(\cT_\fppf) \ra  D(\AAb).
}
In particular, $\bfG_{m,\cT}(U)=\Gamma(U,\cO_U^\tm)$ is an ordinary
abelian sheaf, and big-site pullback gives
\eqn{f^*\bfG_{m,\cT}\cong\bfG_{m,\cT'}.}
We compute units, Picard groups and Brauer groups in these categories.
The ordinary $Rf_*$ need not preserve $D_\cart$; for Cartesian
coefficients we use the adjunctions of Subsection~\ref{cohdes}.

Restriction to the lisse-\'etale site gives, for every abelian sheaf
$\sF$, the cohomology comparison \cite[0GR0]{sp}
\eqn{
 H^i(\cT_\et,\sF)\cong
 H^i(\cT_\liset,\sF|_{\cT_\liset}),\qquad i\ge0.
}
In particular, it applies to $\bfG_m$.
Then for
  $\tau\in\{\et, \fppf\}$,
 we have the cohomology functors
\eqn{
H^i(\cT_\tau, -): \SShAb(\cT_\tau)\ra \AAb.
}
We will write $H_\tau^i(\cT, -)=H^i(\cT_\tau, -)$.
For $X\in\SSch/k$,  $H_\tau^i(X, -)$  recovers the $\tau$-cohomology for
 schemes, cf.   \cite[075F]{sp}.

 Fix  $\sG\in\SShAb((\SSch/k)_\tau)$ and let $\sG_\cT$ be the  pullback of
 $\sG$ to $\cT$.
Then the functor
 $H^i(-_\tau, \sG_-): (\SStk/k)^\op\ra \SSet$ is stable
 (cf. \cite[2.29]{lv2desc}).
If there is  no other confusion, we simply write
 $\sG\in\SShAb(\cT_\tau)$ for $\sG_\cT$.
Thus we can take $F=H_\tau^i(-, \sG)$ in \ref{pt+ob} and
 obtain $\cXAk^{H_\tau^i(\cX, \sG)}$, the $H_\tau^i(-, \sG)$-set.

\egu{[The {\BM}  obstruction]  \label{eg_bm}
The  cohomological  {Brauer--\Grot}  group of an
  algebraic
  stack $\cT$
 is also defined to be
 the {\etale} cohomology group $\Br\cT:=H_\et^2(\cT, \bfG_{m, \cT})$, where
  $\bfG_{m, \cT}$ is the pullback  sheaf of $\bfG_{m, k}$ to $\cT$,
  represented by  $\bfG_{m, k} \tm_k \cT$.
When $\cT$ is a scheme, this agrees with  the cohomological
   {Brauer--\Grot}   group.
The Brauer group of a stack is considered, for example, by Bertolin and
  Galluzzi \cite[Theorem~3.4]{bg19gerbes},  Zahnd \cite[4.3]{zahnd03gerbe-br}
   and Santens \cite{santens2022bm-stack}.
Then, for $F=\Br=H_\et^2(-, \bfG_m)$, we obtain the \emph{{\BM}  set}
  $\cXAk^{\Br}$. If $\cX$ is represented by a scheme,
  then it coincides  with the usual {\BM}  set.
}

\subsection{Torsors over algebraic stacks} \label{torsor}
Let $S$ be a scheme.
 Let $\cX\in \SStk/S$, and $\sG$ be a sheaf of groups
 on $\cX_\tau$, where $\tau\in\{\et, \fppf\}$.
By \cite[Exercise~2.H]{olsson16stack}, we may assume $\cX$ is an object of
 $\cX_\tau$, and  hence $\cX_\tau$ has a final object.

\defi{ \label{defi_torsor}
A (left)
  \emph{$\sG$-torsor over $\cX_\tau$} is a sheaf $\sY$ on $\cX_\tau$ acted on
 by $\sG$,  such that
\enmt{[\upshape (i)]
\item $\sY$ is locally nonempty, i.e., there exists
 a covering $\{\cU_i\ra \cX\}$ in $\cX_\tau$ such that
 $\sY(\cU_i)\neq\emptyset$ for every $i$, and
\item the map
\aln{
  \sG\tm \sY &\ra \sY\tm \sY \\
(g, y) &\mpt (gy, y)
}
 is an isomorphism.

Right  torsors are defined in a similar way.
Unless otherwise stated,  a torsor always means a left torsor.
A left torsor becomes a right torsor by
\eqn{y\cdot g:=g^{-1}\cdot y.}

We denote this  by $\sY\xra{\sG} \cX_\tau$.
}
A morphism between $\sG$-torsors over $\cX_\tau$ is a $\sG$-equivariant
 morphism of sheaves, which is necessarily an isomorphism.

We call $\sY$ \emph{trivial} if it has a global section (or equivalently,
 $\sY\cong \sG$ as $\sG$-torsors over $\cX_\tau$).

Let $\Tors(\cX_\tau, \sG)$ denote the groupoid of all $\sG$-torsors over $\cX_\tau$.
The isomorphism classes of  torsors make
 $\Tors(\cX_\tau, \sG)_{/\cong}$ a pointed set whose  neutral element
 corresponds to the class of trivial torsors.
}

Let  $\check H_\tau^1(\cX, \sG)$ be the {\Cech} cohomology of $\cX_\tau$
 with coefficient $\sG$.
Then  there is a one-to-one correspondence
 (cf. Giraud \cite[III.3.6.5 (5)]{giraud71cohnonab}) of pointed sets
\eq{ \label{eq_tors=ceckH1}
\Tors(\cX_\tau, \sG)_{/\cong} \os{\sim}{\ra} \check H_\tau^1(\cX, \sG).
}
In the case of $\sG\in\SShAb(\cX_\tau)$, we have a bijection (cf.
 \cite[IV.3.4.2 (i)]{giraud71cohnonab})
\eq{ \label{eq_cechH1=H1}
\check H_\tau^1(\cX, \sG) \os{\sim}{\ra} H_\tau^1(\cX, \sG),
}
 which maps the neutral element to $0$ and is functorial both in  $\cX$
 and $\sG$.

Since in particular a sheaf $\sY$ on $\cX_\tau$ is a stack
 (in groupoids) over $\cX_\tau$,
 one obtains a morphism of stacks $\sY\ra \cX$ over $S$ (cf.
  \cite[09WX]{sp}).
If $f: \cZ\ra \cX$ is a morphism of  algebraic $S$-stacks and
 $\sY\in \Tors(\cX_\tau, \sG)$,
 then obviously the pullback
 $\sY_\cZ=\cZ\tm_\cX \sY$ is  also a $\sG$-torsor.
It is the pullback sheaf of $\sY$ under $f^*: \SSh(\cX_\tau)\ra\SSh(\cZ_\tau)$.
This operation corresponds to the cohomological pullback
\eqn{
f^*: \check H_\tau^1(\cX, \sG)\ra \check  H_\tau^1(\cZ, \sG).
}
In particular, if $\sG = G$ where $G$ is an affine $S$-group scheme,
 (viewed as a sheaf of groups on $\cX_\tau$
 by pullback),
 then for any $S$-scheme $T$, since
    $\sY_T=T\tm_\cX \sY$ is   a $G$-torsor in $\tau$ topology,
    $\sY_T$ is a scheme.
This shows that $\sY\ra \cX$ is representable.

\lemm{ \label{lemm_tors_chp}
Let $G$ be an $S$-group scheme.
Then any $\sY\in\Tors(\cX_\tau, G)$ is an object in $\SStk/S$.
}
\pf{
We imitate the proof from  \cite[Lemma~30.8]{olsson07stack-notes}.
To verify the representability of the diagonal $\sY\ra \sY\tm_S \sY$,
 we first show the representability of $\sY\ra \sY\tm_\cX \sY$.
Let  $p_1\tm p_2: T\ra \sY\tm_\cX \sY$ be a $1$-morphism from an $S$-scheme.
One checks that the fiber product $\cZ=\sY\tm_{\sY\tm_\cX \sY} T$ fits into the
 $2$-{\Cart} diagram of $S$-stacks
\eqn{\xymatrix{
\cZ\ar[r]\ar[d] &T\ar[d]^-{(p_1\tm \id_T)\tm(p_2\tm \id_T)} \\
\sY\tm_\cX T\ar[r]^-{\Delta} &(\sY\tm_\cX T)\tm_\cX (\sY\tm_\cX T).
}}
Note that $\sY\tm_\cX T\ra T$ has a section (for example, $p_1\tm \id_T$)
 as a $G$-torsor over $T_\tau$.
Thus it is a trivial torsor in $\Tors(T_\tau, G)$ and in particular
 $\sY\tm_\cX T\cong G\tm_S T$
 is a scheme.
It follows that $\cZ$ is an object in $\SSpc/S$
 (the category of algebraic $S$-spaces,
 in the sense of \cite[025Y]{sp}) and  $\sY\ra \sY\tm_\cX \sY$
 is representable.
Also $\sY\tm_\cX \sY\ra \sY\tm_S \sY$ is representable since it is a base
 change of the diagonal  $\cX\ra \cX \tm_S \cX$.
It follows that the diagonal $\sY\ra \sY\tm_S \sY$ is representable
 since it is the composition
 $\sY\ra \sY\tm_\cX \sY\ra \sY\tm_S \sY$.

To obtain an atlas of $\sY$, one notes that $\sY$ is locally trivial
 (cf. \cite[1.4.1.1]{giraud71cohnonab}).
Then there is a covering $\{U_i\ra \cX\}$ in $\cX_\tau$
 (where we may assume $U_i$ are schemes)  such that
 $\sY\tm_\cX U_i\cong G \tm_S U_i$.
Thus $Y= \coprod_i G \tm_S U_i$ is a scheme and
 $Y\cong \coprod_i \sY\tm_\cX U_i\ra \sY$
 is an \emph{fppf-atlas} for $\sY$.
It follows by
 Artin’s representability theorem
 (cf., for example, \cite[06DC]{sp}) that $\sY\in\SStk/S$.
}

In the following, if not otherwise explained, a
 \emph{$\sG$-torsor over $\cX$} means a
 $\sG$-torsor over $\cX_\fppf$.

Thus, for a flat and locally finitely presented $S$-group scheme $G$
 acting on an algebraic $S$-stack $\cY$ over $\cX$, the morphism
 $\cY\ra \cX$ is a $G$-torsor  if and only if
 it is fppf and
 $(g, y) \mpt (gy, y)$ induces an equivalence
 $G\tm_S\cY\xra{\sim}\cY\tm_\cX\cY$
 \cite[04UV]{sp}.

\rk{ \label{rk_tors_quo}
Let $G$ be a smooth $S$-group scheme and $\cX$
 an algebraic $S$-stack.
If an atlas  $Y\ra \cX$  is a $G$-torsor, then we have an equivalence
 $\cX\xra{\sim} [Y/G]$ of algebraic $S$-stacks.
Conversely, if $Y$ is an algebraic $S$-space acted on by $G$,
 and we have $\cX\xra{\sim} [Y/G]$, then $Y\ra \cX$ is a $G$-torsor.

Indeed, this is the special case of
 \cite[Exercises  4.4.11 and 4.4.18]{moduli}.
}

\lemm{ \label{lemm_quo_tor}
Let $G\hra H$ be a closed normal subgroup of smooth $S$-group schemes
 such that $F=H/G$ is a smooth $S$-group scheme,
 and let $X$ be an $S$-scheme acted on by $H$.
Then we obtain a $1$-morphism
\eq{ \label{eq_qt_stru}
  f: [X/G]\ra [X/H]
}
 of algebraic
 $S$-stacks, an action
\eq{ \label{eq_qt_act}
 [H/G]\tm_S [X/G]\ra [X/G],
}
 and a natural equivalence of algebraic $S$-stacks
\eq{ \label{eq_qt_triv}
  [H/G]\tm_S [X/G]\xra{\sim}    [X/G]\tm_{[X/H]} [X/G]
}
 defined by $(h, x)\mpt (hx, x)$.

In particular, if we have an exact sequence of
  smooth $S$-group schemes
\eq{ \label{eq_qt_ses}
  1\ra G\ra H\ra F\ra 1,
}
 then $F\cong [H/G]$ and  the action \eqref{eq_qt_act}
 makes \eqref{eq_qt_stru} an $F$-torsor.
}
\pf{
Since we have  a homomorphism  $G\ra H$,
 the $1$-morphism \eqref{eq_qt_stru}
 and the action \eqref{eq_qt_act}
 follow from functoriality
 of forming quotient stacks.

To show the equivalence \eqref{eq_qt_triv},
 let $f_H : X\ra [X/H]$
 and $f_G : X\ra [X/G]$  be the natural $1$-morphisms (i.e.,
  the presentations of  quotient stacks).
Then we have $f\circ f_G = f_H$.
Consider the $2$-commutative diagram of algebraic $S$-stacks
\eqn{ \xymatrix{
  H\tm_S X  \ar@/_2cm/[dd] \ar[r]^-{(h,x)\mpt (hx, x)}
   \ar@{.>}[d]^-g &X\tm_S X
   \ar[d] _-{f_G\tm f_G}
    \ar@/^2cm/[dd]^-{f_H\tm f_H} \\
  [X/G]\tm_{[X/H]} [X/G] \ar[r]\ar[d] & [X/G]\tm_S [X/G]
    \ar[d]_-{f\tm f} \\
  [X/H]\ar[r]^-{\Delta} & [X/H]\tm_S [X/H],
}}
 where the lower square is clearly $2$-{\Cart}, and
 the outer one is $2$-{\Cart} by \cite[Exercise 4.4.11]{moduli}.
It follows that, up to equivalence,
 there is a unique $1$-morphism (shown by the dotted arrow)
 rendering the whole diagram   $2$-commutative and making the
 upper square also $2$-{\Cart}.
Now by Remark \ref{rk_tors_quo}, $f_G\tm f_G$ is a $G\tm_S G$-torsor.
It follows that $g: H\tm_S X\ra [X/G]\tm_{[X/H]}[X/G]$
 is also a $G\tm_S G$-torsor.
Again using Remark \ref{rk_tors_quo},  we obtain that
\eqn{
  [X/G]\tm_{[X/H]}[X/G]\xra{\sim} [(H\tm_S X)/(G\tm_S G)]\xra{\sim}
   [H/G]\tm_S [X/G],
}
 as desired.

If we have the exact sequence \eqref{eq_qt_ses},
 then $H\ra F$ is a $G$-torsor, and $[H/G]\cong F$ is a
 smooth $S$-group scheme.
Also $f$ (i.e. \eqref{eq_qt_stru})  is fppf.
By \eqref{eq_qt_triv}, $f$ is an $F$-torsor.
The proof is complete.
}

For an algebraic stack $\cX$ over $k$ and a linear $k$-group $G$,
 \emph{twists}
 of $G$-torsors  over $\cX$
 can also be defined via contracted product
 (cf. \cite[3.24 (resp. 3.30)]{lv2desc}).

Let $f: \cY\xra{G} \cX$ be a torsor.
Since  the functor
 $\check H_\fppf^1(-, G): (\SStk/k)^\op\ra \SSet$ is also  stable by
 \cite[3.14]{lv2desc},
 we may define
  $\cXAk^{\check H_\fppf^1(-, G)}$.
One  may also define $\cXAk^f\subseteq \cXAk$ to be
 the obstruction set given by $[\cY]$ (cf. \eqref{eq_XAkA}).
It contains $\cXk$ and satisfies
\eq{ \label{eq_desc}
\cXAk^f=\bigcup_{\s\in H^1(k, G)} f_\s(\cY_\s(\bfA_k)),
}
 where $f_\s: \cY_\s\ra \cX\in\Tors(\cX, G_\s)$ is the twist of
 $\cY$ by $\s$.
See \cite[3.20]{lv2desc} for a proof in the $\SStk/k$ case,
 which is similar to that in \cite[Theorem~8.4.1]{poonen17rational}.

If $\cX$ is a $k$-scheme, the definitions and facts of this
 subsection coincide with classical ones.

The following statements are the field cases of
 \cite[Corollary~3.9 and Theorem~3.8]{lwz26comparison}.
\cor{ \label{cor_contracted_rep}
Let $Y$ be a smooth quasi-projective geometrically integral $k$-variety,
 $G$ a linear $k$-group acting on $Y$, and $\cX=[Y/G]$.
\enmt{[\upshape (i)]
\item \label{it_twist_rep} For $\s\in H^1(k,G)$, let $P_\s$ be the
 corresponding right $G$-torsor over $k$.
Then
\eqn{
 Y_\s=P_\s\tm_k^G Y=(P_\s\tm_kY)/G,
 \qquad a\cdot(p,y)=(pa^{-1},ay),
}
 is a smooth quasi-projective geometrically integral $k$-variety, and
 $f_\s:Y_\s\xra{G_\s}\cX$ is the twisted torsor.
\item \label{it_H_rep} Let $\iota:G\hra H$ be a closed immersion into
 a linear $k$-group. Then
\eqn{
 Z_H:=H\tm_k^G Y=(H\tm_kY)/G,
 \qquad a\cdot(h,y)=(h\iota(a)^{-1},ay),
}
 is a smooth quasi-projective $k$-variety, and
 $g_H:Z_H\xra{H}\cX$ is an $H$-torsor.
If $H$ is connected, then $Z_H$ is geometrically integral.
}}

\subsection{Cohomological descent} \label{cohdes}

Let $\La$ be a ring.
Recall that Liu-Zheng \cite[(5.7)]{lz24enhanced} constructed
 a functor that associates to each $X\in\Chp/S$ a
 closed symmetric monoidal presentable stable $\infty$-category
 $\sD(X,\La)$, and to each $f:Y\ra X$ a functor
 $f^*:\sD(X,\La)\ra\sD(Y,\La)$ with a right adjoint $f_*$.
By \cite[Corollary~5.3.8 and Remark~5.3.10]{lz24enhanced}, its homotopy category is
\eqn{
 h\sD(X,\La)\simeq D_\cart(X_\liset,\La)
 \simeq D_\cart(X_\et,\La).
}
On homotopy categories we write the adjunction as
\eq{ \label{eq_liset}
f^*:D_\cart(X_\liset,\La)\rla D_\cart(Y_\liset,\La):Rf_*.
}
Here $Rf_*$ is the right adjoint in the Cartesian categories; in
\eqref{eq_biget} it is the ordinary derived direct image.
We write $f_*$ for the right adjoint in $\sD$.

\rk{ \label{rk_liset=biget}
The big-\'etale description in \cite[Remark~5.3.10]{lz24enhanced} identifies
\eqref{eq_liset} with
\eqn{
f^*:D_\cart(X_\et,\La)\rla D_\cart(Y_\et,\La):Rf_*.
}
For smooth representable morphisms, inverse image agrees with ordinary
lisse-\'etale inverse image \cite[Corollary~5.3.8]{lz24enhanced}.
For $\bfG_m$ we use the ordinary cohomology comparison
\cite[0GR0]{sp}, which holds for every abelian sheaf.
On bounded-below Cartesian complexes, ordinary and Liu-Zheng cohomology
agree by smooth descent
\cite[06XJ]{sp}, \cite[Theorem~6.2.13(1)]{lz24enhanced}.
}

Let $\bDel$ (resp. $\bDel_+$)
 be the  \emph{category of combinatorial simplices}
 whose objects are linearly ordered sets $[n] = \{0, 1, \dots, n\}$
 for $n\ge 0$ (resp. $n\ge -1$).
Let $f: X_0\ra X_{-1}$ be
 a $1$-morphism in $\Chp/S$, and
 $X_\bu: \bDel_+^\op \ra \Chp/S$ its \emph{{\Cech} nerve},
 namely, a collection of $\{X_n\in \Chp/S\}_{n\ge-1}$
 with simplicial structure and
 $X_n$ being $(n+1)$-fold $2$-fibred product
 $X_0\tm_{X_{-1}} X_0\tm_{X_{-1}}\cdots \tm_{X_{-1}} X_0$.
Denote by $f_n$ the obvious $1$-morphism $X_n\ra X_{-1}$.

\prop{ \label{prop_cohdes}
Let $f: X_0\ra X_{-1}$ be
 a smooth surjective $1$-morphism in $\Chp/S$, and
 $X_\bu: \bDel_+^\op \ra \Chp/S$ its {\Cech} nerve.
\enmt{[\upshape (i)]
\item \label{it_D_lim}
The natural map
\eqn{
\sD(X_{-1}, \La)\ra \varprojlim_{n\in\bDel}\sD(X_n, \La)
}
 is an equivalence of $\infty$-categories,
 where the transition maps in the limit are provided by $*$-pullback.

\item \label{it_K_lim}
For every $K\in \sD(X_{-1}, \La)$, we have a
  cosimplicial object
 $\bDel \ra \sD(X_{-1}, \La)$ whose value at
 $[n]$ is ${f_n}_*f_n^*K$ such that $K$ is its limit.

\item \label{it_E1}
If $f$ is representable by algebraic spaces,
 then for every $\sF\in \SShAb({X_{-1}}_\et)$,
 there is a convergent spectral sequence
\eqn{
  E^{p,q}_1=H^q_\et(X_p, \sF)\Ra H_\et^{p+q}(X_{-1}, \sF),
}
 whose $E_2$-page is
\eq{ \label{eq_coh_des}
E^{p,q}_2=\check H^p(X_0/X_{-1}, \sH_\et^q(\sF))\Ra H_\et^{p+q}(X_{-1}, \sF).
}
}}

\pf{
  The statement \eqref{it_D_lim} is \cite[Theorem~6.2.13(1)]{lz24enhanced}.
  The statement \eqref{it_K_lim} is a direct consequence of \eqref{it_D_lim}.
  Indeed,  by  \eqref{it_D_lim}, for any $L\in  \sD(X_{-1}, \La)$,
   the $*$-pullback gives an equivalence of mapping spaces
   \eqn{
     \Map_{\sD(X_{-1}, \La)}(L, K)\ra \varprojlim_{n\in\bDel}
      \Map_{\sD(X_n, \La)}(f_n^*L, f_n^*K).
   }
  On the other hand, for any $n\ge 0$, adjunction gives
  \eqn{
    \Map_{\sD(X_n, \La)}(f_n^*L, f_n^*K) =
     \Map_{\sD(X_{-1}, \La)}(L, {f_n}_*f_n^*K).
  }
  It follows that
  \eqn{
     \Map_{\sD(X_{-1}, \La)}(L, K)\xra{\sim} \varprojlim_{n\in\bDel}
     \Map_{\sD(X_{-1}, \La)}(L, {f_n}_*f_n^*K)  \xra{\sim}
     \Map_{\sD(X_{-1}, \La)}(L,  \varprojlim_{n\in\bDel} {f_n}_*f_n^*K).
  }
  This shows  \eqref{it_K_lim}.
  The statement \eqref{it_E1} is (1) of \cite[06XJ]{sp}.
}

\cor{ \label{cor_et_fppf}
Let $S$ be a locally {\Noet} scheme and $G$ be a smooth
 commutative group scheme over $S$.
Let $X$ be an algebraic $S$-stack locally of finite type.
Then for all $i\ge0$, the canonical  map  $H^i_\et(X, G)\ra H^i_\fppf(X, G)$
 is an isomorphism.
}
\pf{
Choose a smooth atlas $U\ra X$ with {\Cech} nerve $U_\bu$.
Proposition~\ref{prop_cohdes}\eqref{it_E1} and its fppf analogue
\cite[06XJ]{sp} give compatible spectral sequences
\eqn{
 E^{p,q}_{1,\tau}=H^q_\tau(U_p,G)
 \Ra H^{p+q}_\tau(X,G),
 \qquad \tau\in\{\et,\fppf\}.
}
The $U_p$ are algebraic spaces, so we are reduced to this case.
Choose an \'etale atlas by a scheme; its nerve consists of schemes,
and the same spectral sequences reduce us to schemes, where the
comparison follows from \cite[III, Th\'eor\`eme~11.7]{grothendieck95brauer}.
}

\defi{
Let $X\in \Chp/S$.
  A sheaf   $\sF\in \SShAb(X_\et, \La)$ is called \emph{constructible}
   if there exists an atlas $f: Y\ra X$ such that $f^*\sF$ is constructible.
}
\cor{ \label{cor_H_fin}
Let $X$ be a quasi-separated algebraic $k$-stack of finite type,
and let $\La$ be a Noetherian ring killed by an integer invertible in $k$.
If $\sF\in \SShAb(X_\et, \La)$ is constructible,
 then  $H^i_\et(\ol X, \sF)$ is a finitely generated $\La$-module
 for every $i\ge 0$.
}
\pf{
Choose a finite-type smooth atlas $U\ra X$ with {\Cech} nerve $U_\bu$.
Proposition~\ref{prop_cohdes}\eqref{it_E1} gives
\eqn{
E_1^{p,q}=H^q_\et(\ol U_p,\sF)
\Ra H^{p+q}_\et(\ol X,\sF),
}
Since $X$ is quasi-separated, each $U_p$ is a quasi-separated algebraic
space of finite type over $k$. Since $\La$ is Noetherian and the
induced filtrations are finite, we are reduced to algebraic spaces.
Choose a finite-type \'etale atlas of such a space by a scheme.
Applying the same spectral sequence reduces the desired result to
schemes, where it follows from \cite[Corollary~9.3.3]{fu11etale}.
}

\section{Calculation of  Brauer groups}
 \label{calc}
Many important algebraic stacks can be realized as quotients of schemes by
 group actions.
In this section, we calculate their Brauer groups.
Since $k$ is a field of characteristic zero, any algebraic $k$-group is smooth.

\subsection{Sansuc's exact sequence} \label{sansuc}
Let $Y$ be a smooth $k$-variety, $G$ be a  connected  linear $k$-group,
 and $X\xra{G} Y$ be a $G$-torsor such that $X$ is geometrically integral.
Sansuc \cite[Proposition~6.10]{sansuc81groupe} gave an exact sequence,
\eqn{
0\ra \UU Y\ra \UU X\ra \hat G(k)\ra \Pic Y\ra
 \Pic X\ra \Pic G\ra \Br Y\ra
 \Br X,
}
 where $\UU=\bfG_m/k^\tm,\ \Pic=H_\et^1(-, \bfG_m) \in \PPSh(\SSch/k)$
  with  $k^\tm$ being  the constant presheaf, and
  $\hat G(k)=\Hom_{k-\tu{groups}}(G,\bfG_m)$ is the group of
  $k$-characters of $G$.
Now we generalize it to quotient stacks.

\thm{ \label{thm_sansuc}
Let $X$ be a smooth geometrically integral
 $k$-variety, and $G$ be a  connected linear $k$-group acting on
 $X$.
Let $\cY=[X/G]$ be the quotient stack, i.e.,
 $f: X\ra\cY$ is a $G$-torsor over an algebraic
 stack $\cY$ $($see Remark  \ref{rk_tors_quo}$)$.
Then we have an exact sequence
\eq{ \label{eq_sansuc}
0\ra \UU\cY\xra{f^*} \UU X\ra \UU G\ra \Pic\cY\xra{f^*}
 \Pic X\ra \Pic G\ra \Br\cY\xra{f^*}
 \Br X\xra{\rho^*-p_2^*} \Br(G\tm_k X),
}
 where $\rho,\ p_2: G\tm_k X\ra X$ are the action of $G$ and
  the projection to $X$, respectively.
}
\pf{
Let $X_\bullet$ be a {\Cech} nerve of $X/\cY$.
For any $\sF\in\SShAb(\cY_\et)$,
 since  $X\ra \cY$ is smooth surjective and representable,
 by  Proposition \ref{prop_cohdes} \eqref{it_E1},
 there is a spectral sequence
\eqn{
E^{p,q}_1=H^q_\et(X_p, \sF)\Ra H_\et^{p+q}(\cY, \sF),
}
 whose $E_2$-page is
\eq{ \label{eq_coh_desc}
E^{p,q}_2=\check H^p(X/\cY, \sH_\et^q(\sF))\Ra H_\et^{p+q}(\cY, \sF).
}
Then we obtain an exact sequence
\eqn{
  0\ra E_2^{1,0}\ra  H_\et^1(\cY, \sF)\ra E_2^{0,1}\ra  E_2^{2,0}\ra
   \ker( H_\et^2(\cY, \sF)\ra E_2^{0,2}) \ra E_2^{1,1}\ra E_2^{3,0}.
}
Taking $\sF=\bfG_m$,
 it becomes
\gan{
  0\ra \check H^1(X/\cY, \bfG_m)\ra \Pic \cY\ra \check H^0(X/\cY, \Pic)\ra
   \check H^2(X/\cY,\bfG_m)\ra \\
   \ker( \Br\cY\ra\check H^0(X/\cY, \Br)) \ra \check H^1(X/\cY, \Pic)\ra
   \check H^3(X/\cY, \bfG_m).
}

Note that for $n\ge 0$, $X_n$ is equivalent to the
 $(n+1)$-fold $2$-fibred product
 $X\tm_\cY X\tm_\cY\cdots \tm_\cY X=G^n\tm_k X$,
 which is a $k$-variety.
Then \cite[Lemme~6.12]{sansuc81groupe}  still applies to
  $F=\UU$ and   $\Pic$, and
 produces the maps $\UU X \ra \UU G$ and $\Pic X \ra \Pic G$.
More precisely,
we obtain that
\eqn{
  \check H^i(X/\cY,\bfG_m) = \check H^i(X/\cY, \Pic) = 0
}
 for $i\ge2$. Thus $E_2^{2,1}=E_2^{3,0}=0$, so
$E_\infty^{0,2}=E_2^{0,2}$ and the edge map
$\Br\cY\ra\check H^0(X/\cY,\Br)$ is surjective.
Together with the preceding exact sequence and
\cite[Lemme~6.12]{sansuc81groupe}, this gives
\gan{
0\ra \check H^1(X/\cY, \bfG_m) \ra \Pic \cY\ra \check H^0(X/\cY, \Pic)\ra 0, \\
0\ra \check H^1(X/\cY, \Pic) \ra \Br \cY\ra \check H^0(X/\cY, \Br)\ra 0, \\
0\ra \check H^0(X/\cY, \UU) \ra \UU X\ra \UU G \ra
 \check H^1(X/\cY, \bfG_m)\ra 0, \\
0\ra \check H^0(X/\cY, \Pic) \ra \Pic X\ra \Pic G \ra
 \check H^1(X/\cY, \Pic)\ra 0,
}
 where $\Br = H_\et^2(-, \bfG_m)\in \PPSh(\SSch/k)$ is viewed as a presheaf.

The rest of the proof is to combine these exact sequences to
 obtain the exact sequence \eqref{eq_sansuc},
 which is exactly the same as the one in the proof of
  \cite[Proposition~6.10]{sansuc81groupe}.
}

\rk{
We have  $U(G)=\hat G(k)$ as in
 \cite[Proposition~6.10]{sansuc81groupe}.
}

\cor{ \label{cor_sansuc_BG}
Let $G$ be a  connected linear $k$-group.
Then we have $\UU BG=0$, $\Pic BG=\UU G$ and
 a split short exact sequence
 $0\ra \Pic G\ra \Br BG\ra \Br k\ra 0$.
}
\pf{
Take $X=\Spec k$ in Theorem  \ref{thm_sansuc} and note that the atlas
 $\Spec k\ra BG$ is a section
 of the structure morphism $BG\ra\Spec k$.
}

In general, for a torsor   $\cX\xra{G}\cY$ where $\cY\in \SStk/k$ and $G$ is a
 $k$-group, let $\rho,\ p_2: G\tm_k \cX\ra \cX$ be
 the action of $G$ and the projection, respectively.
Let $p_1: G\tm_k \cX\ra G$ be the projection to the factor $G$.
Following  \cite[D\'efinition~3.1]{cao18approximation}, we may define the
 \emph{invariant Brauer subgroup} to be
\eqn{
\Br_G\cX\colon = \{ b\in \Br \cX\mid  \rho^*b-p_2^*b\in p_1^*\Br G\}
}
We also recall the following notation
\gan{
\Br_1(-)=\ker(\Br(-)\ra \Br(-\tm_k\ol k)), \\
\Br_a(-)=\coker(\Br k\ra \Br_1(-)), \\
\Br_e(G)=\ker(\Br_1G\xra{e_G^*} \Br k).
}
Note that we have $\Br_e G\xra{\sim} \Br_a G$.

The next statement extends  \cite[Th\'eor\`eme~3.10]{cao18approximation}.
\cor{ \label{cor_sansuc_modi}
With the notation of Theorem \ref{thm_sansuc}, we have an exact subsequence
 of \eqref{eq_sansuc}
\eqn{
\xymatrix{
\cdots\ar[r] &\Pic X\ar[r]\ar@{=}[d] &\Pic G\ar[r]\ar@{=}[d]
 &\Br\cY\ar[r]
 \ar@{=}[d] &\Br_G X\ar[r]\ar@{^(->}[d]\ar@{}[dr]|-{\square}
 &\Br_e(G)\ar@{^(->}[d]^-{p_1^*} \\
\cdots\ar[r] &\Pic X\ar[r] &\Pic G\ar[r] &\Br\cY\ar[r]^-{f^*}
 &\Br X\ar[r]^-{\rho^*-p_2^*} &\Br(G\tm_k X),
}}
 where the right square is {\Cart}.
}
\pf{
%
Since $\cY$ is the quotient stack $X/G$, we have the $2$-{\Cart}
 diagram
\eqn{\xymatrix{
G\tm_k X\ar[r]^-{\rho}\ar[d]^-{p_2} &X\ar[d]^-{f} \\
 X\ar[r]^-{f} &\cY.
}}
In particular, $f\rho\simeq fp_2$, and hence $(\rho^*-p_2^*)f^* = 0$.
It follows that $f^*: \Br \cY\ra \Br X$ has image in $\Br_G X$.

Since $X$ and $G$ are varieties, by \cite[Lemme~6.6]{sansuc81groupe},
 we have $\Br_a(G\tm_k X)\cong \Br_a G \oplus
 \Br_a(X)$ and $\Br_a G\cong \Br_e G$.
It follows that $p_1^*: \Br_e G\ra \Br_a(G\tm_k X)$ is injective.
Also, \cite[Proposition~3.7]{cao18approximation} holds.
This produces the upper morphism of the right square.

The exactness at $\Br_G X$ is due to the one at $\Br X$,
 and
 the right square is {\Cart} by the definition of $\Br_G X$.
The proof is complete.
}

\subsection{Torsionness of the Brauer group} \label{torsion}

\lemm{[Restriction-Corestriction] \label{lemm_res_cores}
Let $f: \cX\ra \cY$ be a representable finite {\etale} $1$-morphism
of constant degree $n$
 between algebraic stacks
 and $\sF\in\SShAb(\cY_\et)$.
Then, for every $i\ge0$, the composition
\eq{ \label{eq_res_cores}
H_\et^i(\cY, \sF)\xra{f^*} H_\et^i(\cX, f^*\sF)\xra{f_*} H_\et^i(\cY, \sF)
}
 is the map of multiplication by $n$.
}
\pf{
The finite {\etale} traces of \cite[03SH]{sp} on scheme base changes
of $f$ are compatible with base change, so they glue to
$\Tr_f:f_*f^*\sF\ra\sF$. Since $f_*$ is exact, the composite
\eqn{
\sF\xra{\eta}Rf_*f^*\sF\xra{\Tr_f}\sF
}
is {\etale} locally
\eqn{
\sF\xra{\Delta}\sF^{\oplus n}\xra{\sum}\sF,
}
hence equals $n\,\id_{\sF}$.
Applying $H^iR\Gamma(\cY_\et,-)$ gives \eqref{eq_res_cores}.
}
\lemm{ \label{lemm_br_torsion_quo}
Let $X$ be a smooth geometrically integral $k$-variety
 and $G$ be a linear $k$-group acting on $X$. Let
 $\cY=[X/G]$ be the quotient stack.
Then the group $\Br\cY$ is torsion.
}
\pf{
Let $G^0$ be the connected component of $G$.
We have the exact sequence of $k$-groups
\eqn{
  1\ra G^0\ra G\ra F\ra 1
}
 where $F$ is finite.
Then by Lemma \ref{lemm_quo_tor}, we have an $F$-torsor
 $[X/G^0]\ra \cY$.
In particular, this is a finite etale $1$-morphism.
By Lemma \ref{lemm_res_cores},
 it suffices to show that $\Br[X/G^0]$ is torsion.
Thus  we may assume $G$ is connected below.

Since the scheme $X$ is smooth integral {\Noet},  $\Br X$ is torsion (cf.
 \cite[II, Proposition~1.4]{grothendieck95brauer}).
Also since $G$ is linear, by \cite[Theorem~1]{rosengarten23picard},
 $\Pic G$ is torsion.
Since $X$ is geometrically integral and $G$ is connected,
Theorem \ref{thm_sansuc} gives the exact sequence
$\Pic G\ra\Br\cY\ra\Br X$.
Thus $\Br\cY$ is torsion.
}

\rk{ \label{rk_nontor}
The linearity assumption in Lemma \ref{lemm_br_torsion_quo} is essential.
For an abelian variety $A$ over $k$ with dual $A^\vee$, we have
\eqn{
\Br BA\cong\Br k\oplus A^\vee(k)
}
by \cite[Lemma~4.2]{shin24br}. Thus $\Br BA$ is not torsion if
$A^\vee(k)$ has positive rank.
}

To obtain a more general torsionness result,
 we need some purity result on Brauer groups. In \cite[Proposition 2.5 iv]{am20brauer}, for regular Noetherian Deligne-Mumford stacks, and in \cite[Theorem 2.1.7]{shin19br}, for regular Noetherian algebraic stacks, there is the following injection on Brauer groups. For convenience of reading, we state the result and  give a proof.
\prop{ \label{prop_MV}
	Let $\cX$ be an algebraic stack over $k$ which can be covered by finitely many
 open substacks $\cX_i=[X_i/G_i]$,
 where $X_i$ is a smooth $k$-scheme and $G_i$ is a
   $k$-group acting on $X_i$ for every $i=1,\dots, n$.
Then we have an injective homomorphism
\eqn{
  \Br \cX \hra \bigoplus_{i=1}^n \Br \cX_i.}}
\pf{
Since $\cX$ is a smooth algebraic
 $k$-stack, by considering a connected component of
 $\cX$,
 we can assume that $\cX$ is a smooth integral stack.
By \cite[0GR0]{sp},
\eqn{
 H^i(\cX_\et,\bfG_m)\cong H^i(\cX_\liset,\bfG_m).
}
All sheaves and direct images in this proof are taken on the ordinary
lisse-\'etale sites. We continue to denote their cohomology by $H^i_\et$;
the displayed comparison identifies the resulting Brauer groups with
those on the big \'etale sites.
Let $\cX\backslash \cX_1$ be the closed substack of $\cX$ with reduced
 substack structure.
Let $(\cX\backslash \cX_1)_{sing}$ be the singular locus of
 the stack $\cX\backslash \cX_1$.
Since $k$ is a field of characteristic zero,
 if $(\cX\backslash \cX_1)_{sing}\neq \emptyset$,
 then it has codimension at least $2$ in $\cX$.
Let $\cU=\cX\backslash (\cX\backslash \cX_1)_{sing}$.
Then $\cX_1=\cX_1\backslash (\cX\backslash \cX_1)_{sing}$ and
 $\cX_1$ is an open substack of $\cU$.
Let $j\colon \cU \to \cX$ and $j_1\colon \cX_1 \to \cU$, be the open morphisms.
By the Leray spectral sequence $H^p_{\et}(\cX, R^qj_*(\bfG_m))
		\Rightarrow H^{p+q}_{\et}(\cU, \bfG_m)$,
 we have an exact sequence
		$$	H_\et^0(\cX, R^1j_*(\bfG_m))\to
     H_\et^2(\cX, j_*(\bfG_m))\to [Ker(H_\et^2(\cU, \bfG_m)\to
     H_\et^0(\cX, R^2j_*(\bfG_m)))].$$
Since  $j\colon \cU \to \cX$ is the open substack morphism,
 by local computation, $R^1j_*(\bfG_m)=0$.
Hence there is an injection
		$$ H_\et^2(\cX, j_*(\bfG_m))\hookrightarrow H_\et^2(\cU, \bfG_m).$$
Similarly, there is an injection
		$$ H_\et^2(\cU, j_{1*}(\bfG_m))\hookrightarrow H_\et^2(\cX_1, \bfG_m).$$
Since $\cX\backslash \cU$ has codimension at least $2$ in $\cX$,
 and $\cX$ is smooth, the natural morphism of sheaves $\bfG_m\to j_*(\bfG_m)$
 is an isomorphism.
Hence there is an injection
 $$H_\et^2(\cX, \bfG_m)\cong
  H_\et^2(\cX, j_*(\bfG_m))\hookrightarrow H_\et^2(\cU, \bfG_m).$$
By local computation, there is an exact sequence of sheaves
 on $\cU_\liset$:
	$$0\to \bfG_m\to j_{1*}(\bfG_m)\to
   \bigoplus_{\cD\in (\cU\backslash \cX_1)^1} i_{\cD*}\ZZ\to 0. $$
Here $(\cU\backslash \cX_1)^1$ is the set of all integral closed substacks of
 $\cU\backslash \cX_1$, which have codimension one in $\cU$.
And $i_{\cD}\colon \cD\to \cU$ is the closed substack morphism.
It induces an exact sequence
		$$ \bigoplus_{\cD\in (\cU\backslash \cX_1)^1} H_\et^1(\cU, i_{\cD*}\ZZ)\to
     H_\et^2(\cU,\bfG_m)\to H_\et^2(\cU, j_{1*}(\bfG_m)).$$
Since $i_{\cD}\colon \cD\to \cU$ is a closed embedding,
 $H_\et^1(\cU, i_{\cD*}\ZZ)\cong H_\et^1(\cD,\ZZ)$.
Let $(\cX\backslash \cX_1)^1$ be the set of all integral closed substacks of
 $\cX\backslash \cX_1$, which have codimension one in $\cX$.
Since $\cX\backslash \cU$ has codimension at least $2$ in $\cX$,
 every $\cD\in (\cU\backslash \cX_1)^1$ is the restriction of one element of
 $(\cX\backslash \cX_1)^1$.
Since $\cU=\cX\backslash (\cX\backslash \cX_1)_{sing},$ $\cD$ is smooth.
Hence $\cD$ is a smooth integral stack.
Similar argument as in \cite[Proposition 2.4.2]{cts21brauer} and \cite[Proposition 5.2]{ls25malle} shows that
  $H_\et^1(\cD,\ZZ)=0.$ It implies that the morphism
		$$H_\et^2(\cU,\bfG_m)\to H_\et^2(\cU, j_{1*}(\bfG_m))$$
  is injective.
Hence the composition of three injections gives an injection
\eqn{
H_\et^2(\cX,\bfG_m)\hookrightarrow H_\et^2(\cX_1, \bfG_m).
}
}

Combining Lemma \ref{lemm_br_torsion_quo} and Proposition~\ref{prop_MV},
 we have the following torsionness result on Brauer groups of stacks.

\cor{ \label{cor_br_torsion}
Let $\cX$ be a regular Noetherian {\DM} stack,
 or an algebraic stack over $k$ which can be covered by finitely many
 open substacks
 $[X_i/G_i]$,
 where $X_i$ is a smooth geometrically integral $k$-variety and $G_i$ is a
  linear $k$-group acting on $X_i$ for every $i$.
Then $\Br \cX$ is torsion.
}
\pf{
The case that $\cX$  is regular Noetherian {\DM} stack follows from
 \cite[Proposition 2.5 (iii)]{am20brauer}. We now show the other case.
By Proposition \ref{prop_MV},
 it suffices to assume  that $\cX=[X/G]$.
Then the result follows from Lemma \ref{lemm_br_torsion_quo}.
}

\section{Descent theory and the {\BM}  pairing} \label{descentfund}
Throughout this section, let  $p:\cX\ra k$ be
 a quasi-separated algebraic stack over a number field $k$ that admits an atlas $X$,
 where $X$ is a geometrically integral $k$-variety.
In particular,   $\cX$ is geometrically integral over $k$.
The following results extend part of the descent theory of
 {\CT} and Sansuc \cite{cs87descente-ii}, as generalized by
 Harari and Skorobogatov \cite[Theorem~8.1]{hs13descent}.

\subsection{Fundamental exact sequence for groups of multiplicative type}
 \label{fund_seq}
Let $\Gal(\ol k/k)\text{-}\MMod$ be the category of discrete
 $\Gal(\ol k/k)$-modules,  which is equivalent to $\SShAb(k_\et)$.
For $p:\cX\ra\Spec k$, let $Rp_*$ denote the ordinary derived
direct image to the small \'etale site $k_\et$. Its adjunction is
\eqn{
p^*=p^{-1}:D(k)\rla D(\cX_\et):Rp_*
}
 where
\eqn{
 D(k)=D(k_\et)=D(\Gal(\ol k/k)\text{-}\MMod).
}
This is obtained by composing the big-site structure morphism with
the comparison to the small \'etale site \cite[0DDK]{sp}.
The coefficient maps are
\eqn{
 p^{-1}\bfG_{m,k}\ra\bfG_{m,\cX},\qquad
 \bfG_{m,k}\ra Rp_*\bfG_{m,\cX}.
}
The second map is adjoint to the first. The first need not be an
isomorphism, and $\bfG_{m,\cX}$ is not in general Cartesian.
Let $S$ be a $k$-group of multiplicative type (which means
$\ol S$ is a subgroup scheme of a product of finitely many copies of $\bfG_m$).
Then its geometric character group
  $\hat S = \Hom_{\ol k-\tu{gp}}(S, \bfG_m)$
  is  a finitely generated $\Gal(\ol k/k)$-module.
Define in $D(k)$
\eqn{
 \KD'(\cX)=\operatorname{Cone}\bigl(
 \bfG_{m,k}[1]\ra(\tau_{\le1}Rp_*\bfG_{m,\cX})[1]\bigr).
}

Let $a\in H^1(k, \hat S)$.
Since we have a canonical isomorphism
 $H^1(k, \hat S)\cong \Ext_k^1(\ZZ, \hat S)$,
 there is  an exact sequence
 $0\ra \hat S\ra M\ra \ZZ\ra 0$ representing $a$.
Now we apply  \cite[Proposition~1.11]{lvcdext}, where we take
\gan{
 \cA = \SShAb(\cX_\et), \quad \cB = \Gal(\ol k/k)\text{-}\MMod, \quad
  \cC = \AAb, \\
 \Psi_1 = \Hom_k(\ZZ, -), \quad \Psi_2 = \Hom_k(M, -), \quad
  \Psi_3 = \Hom_k(\hat S, -), \\
 A = \bfG_{m,\cX}.
}
Take  $\Phi = p_*$
  so that $\Hom_k(-, p_*-) \cong \Hom_{\cX_\et}(p^*-,-)$ and
  there is an adjoint pair
\eqn{
  p^*:D(k)\rla D(\cX_\et):Rp_*.
}
Also take $f$ to be the coefficient map
$\bfG_{m,k}\ra\tau_{\le1}Rp_*\bfG_{m,\cX}$, so that
 $\Delta=\KD'(\cX)$.
In the diagram below, $\bfG_m$ denotes $\bfG_{m,k}$ in $\Ext_k$
and $\bfG_{m,\cX}$ in $\Ext_{\cX_\et}$.
The maps denoted $p^*$ are induced by inverse image and the coefficient map.
It follows by   \cite[Proposition~1.11]{lvcdext}
 that we have the following commutative diagram with exact rows
{\small\eq{ \label{eq_cdextspc}
\xymatrix@-5mm{
 &\Ext_k^1(\hat S, \bfG_m) \ar[r]^-{p^*}\ar[d]_-{a\cup-}
 &\Ext_{\cX_\et}^1(p^*\hat S, \bfG_m) \ar[r]^-\chi\ar[d]_-{p^*(a)\cup-}
 &\Hom_{D(k)}(\hat S, \KD'(\cX))\ar[r]^-\partial\ar[d]_-{a\cup-}
 &\Ext_k^2(\hat S, \bfG_m)\ar[r]^-{p^*}  &\Ext_{\cX_\et}^2(p^*\hat S, \bfG_m) \\
 &\Br k \ar[r] &\Br_1 \cX \ar[r]^-r
 &H^1(k, \KD'(\cX)).
}}}

\thm{ \label{thm_fund_seq}
We have the \emph{fundamental exact sequence}
\eq{ \label{eq_fund_seq}
H^1(k, S)\xra{p^*} H_\fppf^1(\cX, S)\xra{\chi}
 \Hom_{D(k)}(\hat S, \KD'(\cX))\xra{\partial}
 H^2(k, S)\xra{p^*} H_\et^2(\cX, S),
}
 where the map $\chi$ is the \emph{extended type}.
}
\pf{
By \eqref{eq_cdextspc}, it suffices to show that
\eqn{
  H^i_\fppf(\cX, S) \cong H^i_\et(\cX, S) \cong \Ext_{\cX_\et}^i(p^*\hat S, \bfG_m).
}
The proof is the same as that of \cite[Theorem~8.1]{hs13descent},
 except that we must show that $H^p_\et(\cX, G)=H^p_\fppf(\cX, G)$ is
 still valid for $\cX$ being an algebraic stack satisfying the conditions at the beginning
 of Section \ref{descentfund},
 and $G$ a
 commutative $k$-group.
This follows from Corollary \ref{cor_et_fppf}.
The same comparison gives
$H^i_\et(\cX,\hat S)\cong H^i_\fppf(\cX,\hat S)$.
}

\rk{
By  \eqref{eq_tors=ceckH1} and
 \eqref{eq_cechH1=H1},  $H_\fppf^1(\cX, S)$ classifies $S$-torsors over $\cX$,
 and the map $\chi$ associates a torsor class $[Y]$ to its extended type.
The first three terms of the fundamental exact sequence show that
 two $S$-torsors over $\cX$
 have the same extended type if and only if they are isomorphic
 up to a twist.
}

\subsection{The {\BM} pairing}
 \label{bmpairing}
Suppose that $\cX$ is a separated {\DM} stack of finite type over $k$.
Then
 there exists a finite set of places $T$,
 such that $\cX$ extends to
 a separated {\DM} \emph{integral model} $\sX$ of $\cX$ over $\sO_{k, T}$ 
 (that is, a $\sO_{k, T}$-stack $\sX$ such that
 $\sX\tm_{\sO_{k,T}} k \xra{\sim} \cX$
  (see, for example, \cite[Proposition~3.18]{lmb20champs})),
 and
 this can be shown by looking at a smooth groupoid presentation of $\cX$
 (cf. \cite[Exercise~4.4.18]{moduli}).
Moreover, by \cite[Proposition~13.0.2]{christensen20top}, we have
\eq{ \label{eq_adelic_explicit}
\cXAk\cong \us{T'\supseteq T}{\colim}
 (\prod_{v\in T'}\sX(k_v)\tm \prod_{v\not\in T'}\sX(\sO_v)),
}
 where the colimit runs over all finite sets of places  $T'$ containing $T$.

For any algebraic $k$-stack $\cX$ and $x:\Spec\bfA_k\ra\cX$,
\cite[Theorem~2.13 and Corollary~2.14]{cesnavicius15poitou} gives
\eqn{
\Br\cX\xra{x^*}H^2_\et(\bfA_k,\bfG_m)
\cong\bigoplus_v\Br k_v\xra{\sum_v\inv_v}\QQ/\ZZ.
}
Thus the \emph{{\BM} pairing}
\aln{
\BMp{-}{-}: \cXAk \tm \Br \cX & \ra \QQ/\ZZ, \\
(x, A) & \mpt \sum_{v\in\Omega_k}\inv_v A(x_v),
}
where $x_v$ is the pullback of $x$ to $k_v$, is well-defined.
By global reciprocity, the {\BM} set $\cXAk^{\Br}$ defined in
Example~\ref{eg_bm} coincides with the set defined by this pairing, i.e.,
\eqn{
\cXAk^{\Br}  =  \{x\in \cXAk \mid \BMp{x}{A} = 0 \text{ for all $A\in \Br\cX$}\}.
}

\defi{ \label{defi_Be}
Define
\eqn{
\Be\cX= \ker(\Br_a\cX\ra \prod_{v\in\Omega_k} \Br_a\cX_v)
}
 where $\cX_v=\cX\tm_k k_v$.
}
For $A \in \Be \cX$ and $(x_v)\in\cXAk$, one  checks that
\eqn{
\BMp{(x_v)}{A} = \sum_{v\in\Omega_k} \inv_v A(x_v)
}
 does not depend on the choice of $(x_v)$ or on the representative
 for $A$ in $\Br\cX$.
Assuming $\cXAk\neq\emptyset$, we obtain  a well-defined map
 $i=\BMp{(x_v)}{-}: \Be\cX\ra \QQ/\ZZ$.
This notation will be used in Subsection \ref{unitor-n}.


\subsection{Descent} \label{descent}
Let $a\in H^1(k, \hat S)$.
Replacing the upper exact sequence in \eqref{eq_cdextspc} by
 the exact sequence \eqref{eq_fund_seq}  (see
 the proof of Theorem \ref{thm_fund_seq}),
 we obtain the diagram with exact lower row
\eq{\label{eq_cd_chi_r}
\xymatrix{
 & H_\fppf^1(\cX, S)\ar[r]^-\chi \ar[d]^-{p^*(a)\cup-}
  &\Hom_{D(k)}(\hat S, \KD'(\cX))\ar[d]^-{a\cup-} \\
\Br k\ar[r]  &\Br_1 \cX\ar[r]^-r &H^1(k, \KD'(\cX)).
}}

Let $f:\cY\ra \cX$ be an $S$-torsor.
Define $\lambda=\chi([\cY])$, and
\eqn{
\Br_\lambda \cX = r^{-1}(\lambda_*(H^1(k, \hat S)))\subseteq \Br_1 \cX,
}
 where $\lambda_*a = a\cup \lambda$.

\prop{ \label{prop_f=br_lambda}
Let  $p:\cX\ra k$ be a quasi-separated algebraic stack of finite type
 over a number field $k$.
Let $f:\cY\ra \cX$ be an $S$-torsor.
Then we have $\cXAk^f=\cXAk^{\Br_\lambda}$.
}
\pf{
We carefully imitate the original proof
 of \cite[8.12]{hs13descent} as below.
First note that since $p$ is quasi-separated and of finite type,
 the {\BM} pairing
\aln{
\BMp{-}{-}: \cXAk \tm \Br \cX & \ra \QQ/\ZZ, \\
((x_v)_v, A) & \mpt \sum_{v\in\Omega_k}\inv_v A(x_v),
}
 can also be defined, and  we have
\eqn{
\cXAk^{\Br}  =  \{x\in \cXAk \mid \BMp{x}{A} = 0 \text{ for all $A\in \Br\cX$}\}
}
 (cf. Section \ref{bmpairing}).

Choose a finite set of places $T$, containing the archimedean places,
and a smooth multiplicative-type model $\sS$ of $S$ over $\sO_{k,T}$.
Put
\eqn{
P^1(k,S)=\prod_{v\in T}H^1(k_v,S)
\tm\prod_{v\notin T}'H^1(k_v,S),
}
where the restricted product is taken with respect to
$H^1(\sO_v,\sS)$ for $v\notin T$.
These groups inject into $H^1(k_v,S)$: this follows from
$H^1(\sO_v,\sS^\circ)=0$ and the finite \'etale component group
$\sS/\sS^\circ$.
Thus \cite[Theorem~2.18]{cesnavicius15poitou} gives
\eqn{
H^1(\bfA_k,S)\cong P^1(k,S).
}
In particular, $[\cY](x)\in P^1(k,S)$ for $x=(x_v)_v\in\cXAk$.
For $a\in H^1(k,\hat S)$, both classes are unramified almost everywhere,
and their local cup products vanish there since $\Br\sO_v=0$.
We therefore have the commutative diagram
\eq{ \label{eq_cup}
\xymatrix{
H^1(k, \hat S)\ar@{}[r]|-{\tm} \ar@{=}[d]
 & H^1_\fppf(\cX, S) \ar[r]^-{p^*(-)\cup-} \ar[d]^-{(x_v)_v}
  &\Br \cX \ar[d]^-{(x_v)_v}  \\
H^1(k, \hat S)\ar@{}[r]|-{\tm}
 & P^1(k,S) \ar[r]^-{\cup}
  &\bigoplus \Br k_v\ar[r]^-{\inv=\sum_v\inv_v}  &\QQ/\ZZ
}}
It follows that for any
 $a\in H^1(k, \hat S)$,
 we have
\eq{ \label{eq_cup_eq}
(p^*a \cup [\cY])(x) = a \cup ([\cY](x)).
}
Note that the bottom pairing of \eqref{eq_cup}
 induces a homomorphism
\eqn{
\theta:P^1(k,S) \ra \Hom(H^1(k, \hat S), \QQ/\ZZ),
}
 which fits  into
 the {\PT} exact sequence
 (\cite[proof of Proposition~8.12]{hs13descent}),
\eqn{
\cdots\ra H^1(k, S)\ra P^1(k,S) \xra{\theta}
 \Hom(H^1(k, \hat S), \QQ/\ZZ)\ra\cdots.
}
Combining the preceding discussion with \eqref{eq_cup_eq},
 we see that
\eq{ \label{eq_Br_side}
\text{  $\BMp{x}{p^*a \cup [\cY]} = 0$ for all
 $a\in H^1(k, \hat S)$}
}
 if and only if
\eqn{
[\cY](x)\in \im (H^1(k, S)\ra P^1(k,S)),
}
 which is the definition of $x\in \cXAk^f$.

On the other hand, \eqref{eq_cd_chi_r} shows that
\eqn{
r^{-1}(\lambda_*a)=p^*a\cup[\cY]
+\im(\Br k\ra\Br\cX).
}
By global reciprocity, \eqref{eq_Br_side} is therefore equivalent to
$\BMp{x}{b}=0$ for all $a\in H^1(k,\hat S)$ and
$b\in r^{-1}(\lambda_*a)$, which is to say that
$x\in\cXAk^{\Br_\lambda}$ by the definition of $\Br_\lambda$.
The proof is complete.
}

\section{Comparison to some other cohomological obstructions}
 \label{comparision}  
 Harari \cite{harari02groupes}
 compared  the {\BM} obstruction of varieties  with
 some other cohomological obstructions such as obstructions given by
 torsors under connected groups or abelian gerbes.
In this section, we  extend them to some classes of algebraic stacks.

\subsection{Descent (resp. second descent) obstruction on algebraic stacks}
 \label{desc+2desc}
Since for a (resp. commutative) $k$ group $G$,
 $\check H_\fppf^1(-, G)$ (resp. $H_\et^2(-, G)$) is  stable,
 the following definitions
 make sense.
\defi{ \label{defi_desc_conn_2desc}
The \emph{descent obstruction} is
\eqn{
\cXAk^\desc=\bigcap_{\text{all linear $k$-groups  $G$}}
 \cXAk^{\check H_\fppf^1(-, G)}.
}
We also define
\eqn{
\cXAk^\conn=\bigcap_{\text{all connected  linear $k$-groups  $G$}}
 \cXAk^{\check H_\fppf^1(-, G)},
}
 and the \emph{second descent obstruction} (cf. \cite[4.2]{lv2desc},
 slightly modified here) is
\eqn{
\cXAk^\sdesc=\bigcap_{\text{all commutative  linear $k$-groups  $G$}}
 \cXAk^{H_\et^2(-, G)}.
}}

\subsection{Comparison of obstructions}
\thm{ \label{thm_Br=2desc}
Let $\cX$ be a
smooth quasi-separated algebraic $k$-stack
of finite type that is either {\DM} or
Zariski-locally the quotient of a smooth geometrically integral
 $k$-variety by a linear $k$-group. Then
\eqn{
\cXAk^{\Br}=\cXAk^{\sdesc}.
}}
\pf{ Since $\bfG_m$ is a commutative group, we have $\cXAk^{\sdesc}\subseteq
 \cXAk^{\Br}.$ We need to prove the other side $\cXAk^{\Br}\subseteq
 \cXAk^{\sdesc}$.

Since every commutative linear $k$-group is a product of groups of
 multiplicative type and some $\bfG_a$ factors,
 and $H^2(k,\bfG_a)=H^2(k_v,\bfG_a)=0$ for all places $v\in \Omega_k,$ we have
\aln{
\cXAk^\sdesc
&=\bigcap_{\text{commutative linear $k$-group $G$}}
\cXAk^{H_\et^2(-,G)}\\
&=\bigcap_{\text{$k$-group $G$ of multiplicative type}}
\cXAk^{H_\et^2(-,G)}.
}

Since a $k$-group of multiplicative type is an extension of a finite abelian group
by a torus, it splits over a finite extension of $k$.
Since $\cX$ is smooth and quasi-separated of finite type, it is
regular and Noetherian \cite[0510]{sp}; these properties persist
after a finite extension of $k$.
For a $k$-group $S$ of multiplicative type, by Corollary \ref{cor_br_torsion}
 and restriction-corestriction argument (cf. Lemma \ref{lemm_res_cores}),
 the group $H^2(\cX,S)$ is torsion.
Furthermore, for any fixed $f\in H^2(\cX,S)$,
 there exists some positive integer $n$ such that $f$ lies in the image of
 $H^2(\cX, _nS)\to  H^2(\cX,S)$,
 here $ _nS\colon =\ker(S\xra{n\cdot} S).$
We have
\aln{
\cXAk^\sdesc
&=\bigcap_{\text{$k$-group $G$ of multiplicative type}}
\cXAk^{H_\et^2(-,G)}\\
&=\bigcap_{\text{finite commutative $k$-group $G$}}
\cXAk^{H_\et^2(-,G)}.
}
We assume that $G$ is a finite abelian group over $k$.
Let $\hat{G}$ be the dual of $G$.

 Consider  $H^0(k,\hat{G})$, a finite constant abelian group over $k.$
Then $H^0(k,\hat{G})\cong \bigoplus\limits_{i=1}^r\ZZ/n_i\ZZ$ with integers $r>0$
 and $n_i\in \NN^*$. Let $F$ be the finite abelian $k$-group such that $\hat{F}=H^0(k,\hat{G})$. Then
 $F\cong \bigoplus\limits_{i=1}^r\mu_{n_i}$ has a canonical inclusion
 $F\cong \bigoplus\limits_{i=1}^r\mu_{n_i}\hookrightarrow
  \bigoplus\limits_{i=1}^r\bfG_m=\bfG_m^r$.
Since $\hat{F}=H^0(k,\hat{G})$ is a subgroup of $\hat{G}$, it induces a map
 $G\cong \hat{\hat{G}} \to \hat{\hat{F}}\cong F$.
We have the following commutative diagram
\eqn{\xymatrix@-2mm{
			H^2(k, G) \ar[d]\ar[r]&	H^2(k, F) \ar[d]\ar[r] &H^2(k,\bfG_m^r)\ar[d]\\
			\bigoplus\limits_{v\in\Omega_k}H^2(k_v, G)\ar[r]\ar[d]^{\sum j_v^G}&
       \bigoplus\limits_{v\in\Omega_k}H^2(k_v, F)\ar[r]\ar[d]^{\sum j_v^F}
			&\bigoplus\limits_{v\in\Omega_k}H^2(k_v,\bfG_m^r)\ar[d]^{(\sum j_v)^r} \\
			H^0(k, \hat{G})^*=\hat{F}^*\ar@{=}[r]\ar[d]&
       H^0(k, \hat{F})^*=\hat{F}^*\ar@{^(->}[r]\ar[d]
			&H^0(k,\ZZ^r)^*=(\QQ/\ZZ)^r\ar[d]\\
			0&	0
			&0.
}}
Here the superscript $*$ on an
abelian group means $\Hom(-,\QQ/\ZZ)$.
By the {\PT} exact sequence (\cite[Chapter~I 4.10]{Mi06}),
 every column of this diagram is exact.
By taking a point $(P_v)\in \cXAk^{\Br}$, we have
 the following commutative diagram
\eqn{\xymatrix@-2mm{
			 		H_\et^2(\cX, G) \ar[d]^{(P_v)}\ar[r]&	H_\et^2(\cX, F) \ar[d]^{(P_v)}\ar[r]
           &H_\et^2(\cX,\bfG_m^r)\ar[d]^{(P_v)}\\
			 		\bigoplus\limits_{v\in\Omega_k}H^2(k_v, G)\ar[r]\ar[d]^{\sum j_v^G}&
           \bigoplus\limits_{v\in\Omega_k}H^2(k_v, F)\ar[r]\ar[d]^{\sum j_v^F}
			 			&\bigoplus\limits_{v\in\Omega_k}H^2(k_v,\bfG_m^r)\ar[d]^{(\sum j_v)^r} \\
			 		H^0(k, \hat{G})^*=\hat{F}^*\ar@{=}[r]&
           H^0(k, \hat{F})^*=\hat{F}^*\ar@{^(->}[r]
			 			&H^0(k,\ZZ^r)^*=(\QQ/\ZZ)^r.
}}
Combining these two diagrams,
 we have
\eqn{
\cXAk^{\Br}\subseteq\cXAk^{H_\et^2(\cX,G)}.
}
Hence $\cXAk^{\Br}=\cXAk^{\sdesc}$.
}

\rk{
The last part of this proof is essentially from \cite[Proposition~1]{harari02groupes}.
For the convenience of reading, we give the complete details.
}

\thm{ \label{thm_Br=conn}
Let $\cX$ be an algebraic $k$-stack
 as in Theorem \ref{thm_Br=2desc}.
Then
	\eqn{
		\cXAk^{\Br}\subseteq \cXAk^{\conn}.
}}
\pf{
In view of Theorem \ref{thm_Br=2desc},
 the proofs of \cite[Corollaire~1 and Section~3]{harari02groupes}
	remain valid here.
The 	result follows.
}

\section{Descent for {\BM} sets along a torsor} \label{descentbm}
In this section we extend the results  of Cao \cite[Section~5]
 {cao18approximation}
 on  descent for {\BM} sets along a torsor
 to allow  the base space $X$ to be a quotient stack.
Let $k$ be a number field.
\lemm{ \label{lemm_descent_ker}
Let $Y$ be a smooth geometrically integral quasi-projective
 $k$-variety and $G$ a  connected linear $k$-group acting on $Y$.
Let $f: Y\ra\cX$ be a $G$-torsor over $\cX$ so that
 $\cX=[Y/G]$,
 and $f^*: \Br\cX\ra \Br Y$ be the pullback.
Then we have
\eqn{
\cXAk^{\ker f^*} = \bigcup_{\s\in H^1(k,G)} f_\s(Y_\s(\bfA_k)).
}}
\pf{
Using Theorem \ref{thm_sansuc}, the argument in
\cite[proof of Theorem~3.5]{cdx19comparing} gives
\eqn{
f_\s^*(\ker f^*)\subseteq\im(\Br k\ra\Br Y_\s).
}
Global reciprocity proves the inclusion $\supseteq$.

For the reverse inclusion, choose a closed immersion $\iota:G\hra\SL_n$ as in
 \cite[Lemma~3.10 and its proof]{lwz26comparison}.
By Corollary \ref{cor_contracted_rep}\eqref{it_H_rep}, put
\eqn{
 Z=\SL_n\tm_k^G Y,
 \qquad g:Z\xra{\SL_n}\cX,
 \qquad \pi:\SL_n\ra W:=\SL_n/G.
}
Then $Z$ is smooth, quasi-projective and geometrically integral.
The map $\pi$ induces the first arrow in
\eqn{
 p:Z\ra\cX\tm_kW\ra W,
 \qquad [s,y]\mpt sG.
}
For $w\in W(k)$, writing $P_w=\pi^{-1}(w)$, we have
\eqn{
 p^{-1}(w)\cong P_w\tm_k^G Y\cong Y_{\partial w},
}
 where $\partial:W(k)\ra H^1(k,G)$ is the connecting map of the right
 $G$-torsor $\pi:\SL_n\ra W$.
It sends $w$ to the class $[P_w]$ of the pullback torsor
$P_w=\Spec k\tm_{w,W}\SL_n$
\cite[I, Section~5.4, Proposition~36]{serre97galois}.
Every element of $H^1(k,G)$ lies in $\partial(W(k))$, since
$H^1(k,\SL_n)=\{*\}$.
Corollary \ref{cor_contracted_rep}\eqref{it_twist_rep} shows that every
$Y_\s$ is a smooth quasi-projective geometrically integral $k$-variety.

Write $A=\ker f^*$. It is obvious by $H^1(\bfA_k,\SL_n)=\{*\}$ that
\eq{ \label{eq_XA=gZA}
\cXAk^A = g(Z(\bfA_k)^{g^*A}).
}
Thus to complete the proof,  we only need to show  that each
 $x\in \cXAk^A$ lifts to $Y_\s(\bfA_k)$  for some $\s\in
 H^1(k, G)$.
For this purpose, we use strong approximation of $W$ with respect to
 the {\BM} obstruction (see for example \cite[Theorem~3.3]{colliot2009brauer}),
  which states that
\eqn{
\SL_n(\bfA_k)\cdot W(k)=W(\bfA_k)^{\Br}. 
}
We claim that
\eqn{
p(Z(\bfA_k)^{g^*A})\subseteq W(\bfA_k)^{\Br}.
}
Then for any  $x\in \cXAk^A$,  by \eqref{eq_XA=gZA}
 there is $z\in Z(\bfA_k)^{g^*A}$ such that $x=g(z)$.
Then by the claim we have
\eqn{
p(z)\in  W(\bfA_k)^{\Br} = \SL_n(\bfA_k)\cdot W(k).
}
It follows that $p(z)=s\cdot w$ with some $s\in\SL_n(\bfA_k)$ and $w\in W(k)$.
Let $z'=s^{-1}z$. Then we know that $x=g(z)=g(z')$ and $p(z')=w$.
Thus there exists $y\in p^{-1}(w)(\bfA_k) \cong Y_{\partial w}(\bfA_k)$
 lifting $z'$ and $f_{\partial w}(y)=x$,
 which is the demanded result.

To justify the claim, we only need two facts:
\enmt{[\upshape (i)]
\item $\Br Z\cong \Br \cX$, which follows from  Theorem  \ref{thm_sansuc} and
 $\Br(\SL_n\tm_k Z)\cong \Br Z$ \cite[Lemma~2.1]{cdx19comparing}, and
\item $\Br W\ra \Br Z\ra \Br Y$ factors through  $\Br k$, which follows from
 $Y\cong Z\tm_{W, 1_W} k$, i.e., $p^{-1}(1_W) \cong Y$.
}
The proof is complete.
}

\lemm{ \label{lemm_descent_sha}
Let $f: Y\ra\cX$ be as in Lemma \ref{lemm_descent_ker}.
Assume that
   the {\TS}   set $\Sha^1(G/k) = \{*\}$.
Suppose $A\subseteq \Br\cX$ is a subgroup containing $\ker f^*$ and
 for each $\s\in H^1(k, G)$,
 we are given $B_\s\subseteq \Br_{G_\s}Y_\s$  such that
  ${f_\s^*}^{-1} B_\s=A$.
Then
\eqn{
\cXAk^A = \bigcup_{\s\in H^1(k,G)} f_\s(Y_\s(\bfA_k)^{B_\s}).
}}
\pf{
The inclusion $\supseteq$ follows from functoriality of the {\BM} pairing.
Lemma \ref{lemm_descent_ker} together with $A\supseteq \ker f^*$ implies that
\eqn{
\cXAk^A = \bigcup_{\s\in H^1(k,G)} f_\s(Y_\s(\bfA_k)^{f_\s^*A}).
}
It suffices to show that $f_\s(Y_\s(\bfA_k)^{f_\s^*A})\subseteq
 f_\s(Y_\s(\bfA_k)^{B_\s})$ (the reverse inclusion is immediate since
 $f_\s^*A\subseteq B_\s$), and one may assume that $\s$ is the neutral  element.
Using Corollary \ref{cor_sansuc_modi} and \cite[Th\'eor\`eme~5.1]{demarche11defaut},
 the remaining argument is the same as \cite[Lemme~5.1]{cao18approximation}.
}

\thm{ \label{thm_descent}
Let $Y$ be a smooth geometrically integral quasi-affine
 $k$-variety, and $G$ be a  connected linear $k$-group.
Let $f: Y\ra\cX$ be a $G$-torsor over $\cX$ so that
 $\cX=[Y/G]$.
Then
\eqn{
\cXAk^{\Br} = \bigcup_{\s\in H^1(k,G)} f_\s(Y_\s(\bfA_k)^{\Br_{G_\s}(Y_\s)}).
}}
\pf{
Identify $R=G/G^u$ with a Levi subgroup of $G$, and put
$q:[Y/R]\ra\cX$ and $f_R:Y\ra[Y/R]$.
Since $G\cong G^u\tm_kR$ as varieties,
homotopy invariance, Theorem~\ref{thm_sansuc} and
\cite[Lemma~2.1]{cdx19comparing} give
\eqn{
q^*:\Br\cX\xra{\sim}\Br[Y/R],
\qquad \Br_GY=\Br_RY.
}
For $F=k,\bfA_k$, the bijection $H^1(F,R)\ra H^1(F,G)$ gives
an $R$-reduction $Q$ of $x^*Y$ for every $x\in\cX(F)$.
The induced $R$-equivariant map $Q\ra Y$ defines a lift of $x$
to $[Y/R](F)$. Thus $q(F)$ is surjective, and the isomorphism $q^*$ implies
\eqn{
q([Y/R](\bfA_k)^{\Br})=\cXAk^{\Br}.
}
For $\rho\in H^1(k,R)$ and its image $\sigma\in H^1(k,G)$,
the contracted products identify
\eqn{
Y_\rho\cong Y_\sigma,\qquad
f_\sigma=q\circ(f_R)_\rho,\qquad
\Br_{R_\rho}Y_\rho=\Br_{G_\sigma}Y_\sigma.
}
Thus it suffices to prove the equality for $f_R$, so we may assume
that $G$ is reductive.

Fix $x\in\cXAk^{\Br}$. By Lemma~\ref{lemm_descent_ker}, choose
$\sigma_0\in H^1(k,G)$ such that
$x\in f_{\sigma_0}(Y_{\sigma_0}(\bfA_k))$.
Twisting permutes the terms in the union. Replacing
$(G,Y,f)$ by $(G_{\sigma_0},Y_{\sigma_0},f_{\sigma_0})$, we may
assume that $x\in f(Y(\bfA_k))$.
By \cite[Proposition~5.6]{cao18approximation},
we find a resolution of $G$ by a connected linear $k$-group $H$
  and a $k$-torus $T$
\eqn{
1\ra G\xra{\psi} H\xra{\phi} T\ra 1
}
 such that $\Sha^1(H/k)=\{1\}$ and $H^3(k,\hat T)=0$.
The pushforward of $Y$ under the map
 $\psi_*: H_\fppf^1(\cX, G)\ra H_\fppf^1(\cX, H)$ is
 an $H$-torsor $g:Z=H\tm_k^G Y\xra{H}\cX$.
The map $\phi$ induces the first arrow in
\eq{ \label{eq_resolution}
 p:Z\ra\cX\tm_kT\ra T,
 \qquad [h,y]\mpt\phi(h).
}
By Corollary~\ref{cor_contracted_rep}, $Z$, $Y_\s$ and $Z_\tau$ are
smooth quasi-projective geometrically integral $k$-varieties for
$\s\in H^1(k,G)$ and $\tau\in H^1(k,H)$.

Since $\Sha^1(H/k)=\{1\}$,
 applying  Lemma \ref{lemm_descent_sha} to $g: Z\xra{H} \cX$ we obtain
\eq{ \label{eq_XABr=fZABr}
\cXAk^{\Br} = \bigcup_{\s\in H^1(k,H)} g_\s(Z_\s(\bfA_k)^{\Br_{H_\s}(Z_\s)}).
}
Also  $H^3(k, \hat T) = 0$ implies that $\Br_a H\ra \Br_a G$ is surjective.
It follows from \cite[Proposition~3.13]{cao18approximation} that
$\Br_HZ\ra\Br_GY$ is also surjective.
By \eqref{eq_XABr=fZABr},
 there is some $\s\in H^1(k,H)$ and  $z\in
 Z_\s(\bfA_k)^{\Br_{H_\s}(Z_\s)}$ such that $x=g_\s(z)$.
Since $x\in f(Y(\bfA_k))\subseteq g(Z(\bfA_k))$ and $x=g_\s(z)$,
the torsors $x_v^*Z$ and $x_v^*Z_\s$ both have $k_v$-points. Hence
\eqn{
\s_v=1\quad(v\in\Omega_k),\qquad
\s\in\Sha^1(H/k)=\{1\}.
}
Thus $z\in Z(\bfA_k)^{\Br_H(Z)}$.
Since $\Br_1T=\Br_TT$ by \cite[Lemme~3.6]{cao18approximation} and $p$ is
equivariant, $p^*(\Br_1T)\subseteq\Br_HZ$. Thus
$p( Z(\bfA_k)^{\Br_{H}(Z)}) \subseteq
 T(\bfA_k)^{\Br_1}$.
Also \cite[Proposition~3.3]{cdx19comparing} gives
\eqn{
 T(\bfA_k)^{\Br_1} = \phi(H(\bfA_k)^{\Br_1})T(k).
}
It follows that $p(z)=\phi(h)t$ with some $h\in H(\bfA_k)^{\Br_1}$ and
 $t\in T(k)$.
Let $z'=h^{-1}z\in
 Z(\bfA_k)^{\Br_H(Z)}$.
Then we know that $x=g(z)=g(z')$ and $p(z')=t$.
Denote by  $\partial: T(k)\ra H^1(k, G)$  the connecting map associated to
 \eqref{eq_resolution}.
Then putting $\tau=\partial t$,
 we know that
 there exists $y\in p^{-1}(t)(\bfA_k) \cong Y_\tau(\bfA_k)$ lifting
 $z'$ and $f_\tau(y)=x$.
On the other hand, \cite[Proposition~3.13]{cao18approximation} gives a surjection
$\Br_HZ\ra\Br_{G_\tau}Y_\tau$ by restriction to
$p^{-1}(t)\cong Y_\tau$, where $G_\tau\cong G$.
Since $z'\in Z(\bfA_k)^{\Br_H(Z)}$, we also have
$y\in Y_\tau(\bfA_k)^{\Br_{G_\tau}Y_\tau}$.
This finishes the proof.
}

\section{{\BM} set under a product} \label{prod}

Let $k$ be a number field.
The  product preservation property of the {\BM} set was first established by
 Skorobogatov and Zarhin \cite{sz14product} for smooth geometrically integral
 projective $k$-varieties, and later by the first author
 \cite{lv20brprodge}
 for open ones.
In this section we extend this result to a class of algebraic stacks over $k$.

\defi{
We define $\SStk'/k\subset \SStk/k$ to be the full sub-$2$-category spanned by
 smooth quasi-separated algebraic $k$-stacks
 of finite type that admit an atlas $X$ which is a smooth geometrically
 integral $k$-variety with $\XA^{\Be}\neq\emptyset$
 $($see Definition \ref{defi_Be}$)$,
 and are either {\DM} or
 Zariski-locally the quotient of a smooth geometrically integral
 $k$-variety by  a linear $k$-group.
}

\thm{ \label{thm_prod}
The canonical map
\eqn{
(\cX\tm_k\cY)(\bfA_k)^{\Br}\xra{\sim} \cXAk^{\Br}\tm \cYAk^{\Br}
}
 is a bijection for any $\cX, \cY\in \SStk'/k$.
}

\subsection{A variant of the {\Kunn} formula} \label{kunneth}
Let $\La$ be a commutative Noetherian ring killed by some positive integer.
All derived operations in this subsection are Liu-Zheng operations
on Cartesian complexes, as in Subsection~\ref{cohdes}.
For $\cX\in\SStk/k$ and $K\in D_\cart(\cX_\et,\La)$, we say that $K$ has
\emph{Tor-amplitude in $[c,d]$} if, for a smooth atlas
$f:U\ra\cX$,
\eqn{
f^*K\otimes_\La^L N\in D^{[c,d]}(U_\et,\La)
\qquad\text{for every }N\in D^0(U_\et,\La).
}
This is independent of the atlas, since the condition is tested on
geometric stalks \cite[08FZ,0DJJ]{sp}.
Finite Tor-amplitude means Tor-amplitude in some finite interval.
For $\cX\in \SStk/k$, let $D_\ctf^b(\cX_\et, \La)$ be the  full
 subcategory of $D_\cart(\cX_\et, \La)$  spanned by
 objects $K$ such that $f^*K\in D_\ctf^b(X_\et, \La)$
 for every atlas $X\ra \cX$.
\lemm{[{\Kunn} formula] \label{lemm_kunneth}
Let $p:\cX\ra k$ and $q:\cY\ra k$ be quasi-separated algebraic
stacks of finite type over a field $k$, and let
$K\in D_\ctf^b(\cX_\et,\La)$ and $L\in D_\ctf^b(\cY_\et,\La)$.
Write $u,v$ for the projections from $\cX\tm_k\cY$ and
$p\tm q=pu=qv$. Suppose that $\La$ is killed by an integer
invertible in $k$. If one of $P=Rp_*K,Rq_*L$ satisfies
\eqn{
P\otimes_\La^L D^{\ge0}(k,\La)\subseteq D^{\ge a}(k,\La)
}
for some integer $a$, then the canonical map
\eqn{
Rp_*K\otimes_\La^L Rq_*L
\ra R(p\tm q)_*(u^*K\otimes_\La^Lv^*L)
}
is an isomorphism.
}
\pf{
For a quasi-separated algebraic space $z:Z\ra k$ of finite type,
$z_*$ preserves small colimits \cite[Proposition~6.3.2]{lz24enhanced}.
The shifts of the perfect permutation sheaves $\La[U]$, with $U/k$
finite \'etale, generate $\sD(k,\La)$ under small colimits
\cite[077J]{sp}. Hence the projection formula \cite[0944]{sp}
holds for every $N\in\sD(k,\La)$. If $F$ has Tor-amplitude $[c,d]$
and $N\in\sD^0(k,\La)$, then
\eqn{
(z_*F)\otimes_\La N\simeq z_*(F\otimes_\La z^*N)
\in\sD^{[c,d+e_Z]}(k,\La),
}
where $e_Z$ is a cohomological-dimension bound for $z_*$
\cite[Lemma~6.3.1]{lz24enhanced}. Thus $z_*F$ has finite Tor-amplitude.

Let $M^\bu:\bDel\ra\sD^{\ge b}(k,\La)$ and suppose
$P\otimes_\La\sD^{\ge0}\subseteq\sD^{\ge a}$.
Passing to the inverse limit of the partial totalizations in
\cite[Remark~1.2.4.3, Proposition~1.2.4.5(4)]{HA}, applied to the opposite
category, gives
\eqn{
F_r:=\operatorname{fib}\bigl(\operatorname{Tot}(M^\bu)
\ra\operatorname{Tot}_{\le r}(M^\bu)\bigr)
\in\sD^{\ge b+r+1},
}
since $\sD^{\ge b+r+1}$ is closed under limits.
Exactness of $P\otimes_\La-$ gives the commutative square
\eqn{
\xymatrix@C=2em{
P\otimes_\La\operatorname{Tot}(M^\bu)
\ar[r]^-\alpha\ar[d]_-{\alpha_r}
&\operatorname{Tot}(P\otimes_\La M^\bu)\ar[d]^-{\beta_r}\\
P\otimes_\La\operatorname{Tot}_{\le r}(M^\bu)
\ar[r]^-{\sim}
&\operatorname{Tot}_{\le r}(P\otimes_\La M^\bu).
}
}
Here the bottom arrow is an equivalence because partial totalizations
are formed from finite limits and retracts. The vertical fibres satisfy
\aln{
\operatorname{fib}(\alpha_r)&\simeq P\otimes_\La F_r
\in\sD^{\ge a+b+r+1},\\
\operatorname{fib}(\beta_r)&\in\sD^{\ge a+b+r+1}.
}
The fibre sequence
\eqn{
\operatorname{fib}(\alpha)\ra\operatorname{fib}(\alpha_r)
\ra\operatorname{fib}(\beta_r)
}
therefore gives
\eqn{
\operatorname{fib}(\alpha)\in
\bigcap_{r\ge0}\sD^{\ge a+b+r+1}=0.
}
Consequently,
\eqn{
P\otimes_\La\varprojlim_{n\in\bDel}M^n
\xra{\sim}\varprojlim_{n\in\bDel}(P\otimes_\La M^n).
}

Choose finite-type smooth atlases with \Cech\ nerves
$f_\bu:X_\bu\ra\cX$ and $g_\bu:Y_\bu\ra\cY$, and set
\aln{
A&=p_*K,& A^m&=(pf_m)_*f_m^*K,\\
B&=q_*L,& B^n&=(qg_n)_*g_n^*L,\\
E&=u^*K\otimes_\La v^*L,& h_{mn}&=f_m\tm g_n.
}
Proposition~\ref{prop_cohdes}\eqref{it_K_lim} and
\cite[Lemma~4.2.2]{lz24enhanced} give
\eq{\label{eq_ku_lim}
\begin{gathered}
A\simeq\varprojlim_{m\in\bDel} A^m,\qquad B\simeq\varprojlim_{n\in\bDel} B^n,\\
(p\tm q)_*E\simeq\varprojlim_{(m,n)\in\bDel\tm\bDel}(p\tm q)_*h_{mn*}h_{mn}^*E.
\end{gathered}
}
Choose a common lower bound $b_0$ for $K,L$. Since pullback is
t-exact and direct image preserves lower bounds,
\eqn{
A^m,B^n\in\sD^{\ge b_0}(k,\La)\qquad(m,n\ge0).
}
Each $B^n$ has finite Tor-amplitude by the first paragraph.
The formula for algebraic spaces gives
\eq{\label{eq_ku_esp}
A^m\otimes_\La B^n\xra{\sim}(p\tm q)_*h_{mn*}h_{mn}^*E.
}
By symmetry, we may assume that $A$ satisfies the hypothesis.
Apply the preceding equivalence first to $(A,B^\bu)$ and then
to $(B^n,A^\bu)$ for each fixed $n$. We obtain
\aln{
p_*K\otimes_\La q_*L
&\simeq\varprojlim_{n\in\bDel}(A\otimes_\La B^n)\\
&\simeq\varprojlim_{(m,n)\in\bDel\tm\bDel}(A^m\otimes_\La B^n)\\
&\simeq\varprojlim_{(m,n)\in\bDel\tm\bDel}(p\tm q)_*h_{mn*}h_{mn}^*E
\simeq(p\tm q)_*E.
}
For algebraic spaces, the same argument uses \'etale scheme atlases
and the scheme formula \cite[Corollary~9.3.5]{fu11etale}; the first
paragraph gives
\eqn{
(z_*F)\otimes_\La\sD^{\ge0}(k,\La)
\subseteq\sD^{\ge c}(k,\La).
}
}

For $\La=\ZZ/n\ZZ$ with $n$ invertible in $k$, put
$C=p_{*,\ZZ}\ZZ_\cX\in\sD^{\ge0}(k,\ZZ)$.
The resolution $\La\simeq{[\ZZ\xra{n}\ZZ]}$ in degrees $-1,0$
and perfect change of coefficients \cite[Proposition~6.2.6]{lz24enhanced}
give, for $M\in\sD^{\ge0}(k,\La)$,
\aln{
p_{*,\La}\La_\cX&\simeq C\otimes_\ZZ^L\La,\\
(p_{*,\La}\La_\cX)\otimes_\La^L M
&\simeq C\otimes_\ZZ^L M\in\sD^{\ge-1}(k,\La).
}
The last inclusion follows on geometric stalks from
$\operatorname{gldim}\ZZ=1$. For the invertible sheaf $T=\La(j)$,
the projection formula \cite[Proposition~6.2.4(2)]{lz24enhanced} gives
\eqn{
p_{*,\La}\La_\cX(j)\simeq(p_{*,\La}\La_\cX)\otimes_\La T.
}
Thus the hypothesis of Lemma~\ref{lemm_kunneth} holds with $a=-1$
for constant coefficients and their Tate twists, including $\mu_n$,
since tensoring with $T$ is t-exact.

\lemm{[Smooth base change] \label{lemm_sbc}
Let
\eqn{
\xymatrix{
\cW\ar[r]^-g\ar[d]_-q &\cZ\ar[d]^-p \\
\cY\ar[r]^-f &\cX
}}
 be a $2$-{\Cart} diagram in $\SStk/k$ where $p$ is smooth.
Then
\eqn{
p^*Rf_*\ra Rg_*q^*:
D_\cart(\cY_\et,\La)\ra D_\cart(\cZ_\et,\La)
}
 is an isomorphism of functors.
}
\pf{
  In view of Remark \ref{rk_liset=biget},
 this is \cite[Corollary~0.1.5]{lz24enhanced}.
}

\subsection{The universal torsor of $n$-torsion} \label{unitor-n}
In this subsection, $Rp_*:D(\cX_\et)\ra D(k)$ is the ordinary
derived direct image of Subsection~\ref{fund_seq}.
For $\mu_n$, ordinary and Liu-Zheng cohomology agree by
Remark~\ref{rk_liset=biget}. We apply this ordinary $Rp_*$ to both
$\mu_n$ and $\bfG_m$ in the diagrams below.
Let $\cX$ be a quasi-separated algebraic $k$-stack of finite type.
Assume that $\cX$ is nonempty and geometrically connected, as is the
case if it admits a geometrically integral $k$-variety as an atlas.
These hypotheses hold for every $\cX\in\SStk'/k$.
By Corollary \ref{cor_H_fin}, where we take $\La = \ZZ/n\ZZ$,
 the $\Gal(\ol k/k)$-module
 $H_\et^1(\ol \cX, \mu_n)$ is finite as an abelian group.
Let $S_\cX$ be the group of multiplicative type dual to  it, that is,
 $\hat S_\cX=H_\et^1(\ol \cX, \mu_n)$.
Then as in \cite[Section~2]{cao20sous},
 we have the following
 commutative diagram of distinguished triangles.
\eqn{
\xymatrix{
\Delta_n[-2] \ar[r]\ar[d] &\mu_n\ar[r]\ar[d] &Rp_*\mu_n\ar[r]^-{+1}\ar[d] & \\
\Delta[-2] \ar[r]\ar[d]^-{\tm n} &\bfG_m\ar[r]\ar[d]^-{\tm n}
 &Rp_*\bfG_m\ar[r]^-{+1}\ar[d]^-{\tm n} & \\
\Delta[-2] \ar[r]\ar[d]^-{+1} &\bfG_m\ar[r]\ar[d]^-{+1}
 &Rp_*\bfG_m\ar[r]^-{+1}\ar[d]^-{+1} & \\
 & &
}}
 where  $\Delta$ (resp. $\Delta_n$) is the cone of
  $\bfG_m[1]\ra Rp_*\bfG_m[1]$
  (resp. $\mu_n[1]\ra Rp_*\mu_n[1]$),
  so that  we have
 $\tau_{\le0} \Delta = \KD'(\cX)$. 
Also we have $H^0(\Delta_n)\cong \hat S_\cX$.
Taking $S=S_\cX$  in Theorem \ref{thm_fund_seq}, we obtain the
 exact sequence
\eq{ \label{eq_fund_seq_S_X}
H^1(k, S_\cX)\xra{p^*} H_\fppf^1(\cX, S_\cX)\xra{\chi}
 \Hom_{D(k)}(\hat S_\cX, \KD'(\cX))\xra{\partial}
 H^2(k, S_\cX)\xra{p^*} H_\et^2(\cX, S_\cX).
}
Truncating the canonical map  $\Delta_n\ra \Delta$,
 we obtain
\eq{ \label{eq_tau}
  \tau_{\le0}(\Delta_n\ra \Delta): \tau_{\le0}\Delta_n
  \ra \tau_{\le0}\Delta.
}
Geometric connectedness gives $\mu_n\xra{\sim}p_*\mu_n$, so the first
triangle gives $H^i(\Delta_n)=0$ for $i<0$.
Hence $\tau_{\le0}\Delta_n\simeq\hat S_\cX$, while
$\tau_{\le0}\Delta = \KD'(\cX)$.
It follows that  \eqref{eq_tau} is an element
 of  $\Hom_{D(k)}(\hat S_\cX, \KD'(\cX))$.

\defi{ \label{defi_univ_n_tor}
A \emph{universal torsor  of $n$-torsion} of $\cX$
 is an $S_\cX$-torsor over $\cX$ whose class in
 $H_\fppf^1(\cX, S_\cX)$ is  in the inverse image
 of \eqref{eq_tau}
 under $\chi$  (see \eqref{eq_fund_seq_S_X}).
}

\lemm{ \label{lemm_type_surj}
Suppose that $\cX/k$ is quasi-separated and admits an atlas $X$ which is a smooth geometrically
 integral $k$-variety and satisfies $\XA^{\Be}\neq\emptyset$.
Let $S$ be a group of multiplicative type.
Then
\eqn{
\chi:  H_\fppf^1(\cX, S)\ra \Hom_{D(k)}(\hat S, \KD'(\cX))
}
 in \eqref{eq_fund_seq} is surjective.
}
\pf{
Recall the exact sequence \eqref{eq_fund_seq}
\eqn{\xymatrix{
H^1(k, S)\xra{p^*} H_\fppf^1(\cX, S)\xra{\chi}
 \Hom_{D(k)}(\hat S, \KD'(\cX))\xra{\partial}
 H^2(k, S)\xra{p^*} H_\et^2(\cX, S).
}}
First note that $\chi$ is surjective if and only if
 $H^2(k, S)\ra H_\fppf^2(\cX, S)$ is injective.
Let $f: X\ra \cX$ be the given atlas such that
 $X$  is a smooth geometrically
 integral $k$-variety with $\XA^{\Be}\neq\emptyset$.
It follows from \cite[Corollary~8.17]{hs13descent} that
 $H^2(k, S)\ra H_\fppf^2(X, S)$ is injective,
 which is the composition
\eqn{
  H^2(k, S)\xra{p^*} H_\fppf^2(\cX, S)\xra{f^*} H_\fppf^2(X, S).
}
It follows that $p^*$ is also injective.
The proof is complete.
}

By Lemma \ref{lemm_type_surj},  if $\cX$ (resp. $\cY$)
 satisfies the assumptions,
 then there exists a universal torsor  of $n$-torsion
 $\cT_\cX\xra{S_\cX}\cX$ of $\cX$
 (resp. $\cT_\cY\xra{S_\cY}\cY$ of $\cY$).
Then the homomorphism introduced  by Skorobogatov and Zarhin (cf. \cite[(1-4)]
 {cao23sous-corr}) can also be defined, that is,
\aln{
\epsilon': \Hom_k(S_\cX\otm_{\ZZ/n\ZZ} S_\cY, \mu_n) &\ra
 H_\et^2(\cX\tm_k \cY, \mu_n) \\
\phi &\mpt \phi_*(\cT_\cX\cup\cT_\cY).
}
\lemm{ \label{lemm_kunneth_H}
Let $\cX$ and $\cY$ be quasi-separated algebraic $k$-stacks satisfying
 the assumptions of Lemma \ref{lemm_type_surj}.
Then we have isomorphisms of $\Gal(\ol k/k)$-modules
\ga{
(p_1^*, p_2^*): H_\et^1(\ol\cX, \mu_n)\oplus H_\et^1(\ol\cY, \mu_n)  \xra{\sim}
 H_\et^1(\ol\cX\tm_{\ol k} \ol\cY, \mu_n), \label{eq_H1} \\
(p_1^*, p_2^*, \epsilon'): H_\et^2(\ol\cX, \mu_n) \oplus
 H_\et^2(\ol\cY, \mu_n)
 \oplus
 \Hom_{\ol k}(S_{\ol\cX}\otm_{\ZZ/n\ZZ} S_{\ol\cY}, \mu_n) \xra{\sim}
 H_\et^2(\ol\cX\tm_{\ol k}\ol\cY, \mu_n). \label{eq_H2}
}}
\pf{
Note that $\cX$ and $\cY$ are automatically smooth over $k$.
The original proof (cf. \cite[Th\'eor\`eme~2.1]{cao23sous-corr}) is also available
 if we use stacky versions of the lemmas
 (Lemmas \ref{lemm_kunneth} and \ref{lemm_sbc}).
}

\subsection{Proof of the product preserving property} \label{prodpf}
Now we can prove Theorem \ref{thm_prod}.
We need  to verify that
\eqn{
(\cX\tm_k\cY)(\bfA_k)^{\Br}\xra{\sim} \cXAk^{\Br}\tm \cYAk^{\Br}
}
 for any $\cX, \cY\in \SStk'/k$.
We go through the proofs in \cite[Section~5]{sz14product}.
By Corollary~\ref{cor_br_torsion}, $\Br\cX$ and $\Br\cY$ are torsion,
as is $\Br(\cX\tm_k\cY)$ when both factors are {\DM} or both
are locally quotients.
In the mixed case, say $\cX$ is {\DM}, choose a dense open
$\cU\subseteq\cX$ and a finite \'etale cover $V\ra\cU$ by a smooth
variety \cite[proof of Proposition~2.5(iii)]{am20brauer}.
Choose a finite extension $k'/k$ such that the components of $V_{k'}$
are geometrically integral. Corollary~\ref{cor_br_torsion} shows that
$\Br(V_{k'}\tm_{k'}\cY_{k'})$ is torsion.
Since
\eqn{
\Br(\cX\tm_k\cY)\hra\Br(\cU\tm_k\cY)
}
by \cite[Theorem~5.3]{ls25malle}, Lemma~\ref{lemm_res_cores} applied to
$V_{k'}\tm_{k'}\cY_{k'}\ra\cU\tm_k\cY$ shows that
$\Br(\cX\tm_k\cY)$ is torsion.
Also, there exist universal torsors  of $n$-torsion
 $\cT_\cX\xra{S_\cX}\cX$ of $\cX$ and
 $\cT_\cY\xra{S_\cY}\cY$ of $\cY$,
 since both $\cX$ and $\cY$ satisfy
 the assumptions of Lemma \ref{lemm_type_surj}.
Then we use \cite[Lemmas 5.2 and 5.3]{sz14product} to finish the proof.
In the proof of \cite[Lemma~5.2]{sz14product}, we replace
\cite[Corollary~1.8 and Theorem~2.6]{sz14product} by
\eqref{eq_H1} and \eqref{eq_H2}, respectively.
The proof is complete.

\cor{ \label{cor_prod_quo}
If $\cX$ and $\cY$ are stacks that are quotients of smooth geometrically integral
 quasi-affine $k$-varieties by connected linear $k$-groups, then
\eqn{
(\cX\tm_k\cY)(\bfA_k)^{\Br}\xra{\sim} \cXAk^{\Br}\tm \cYAk^{\Br}
}}
\pf{
Combine Theorems \ref{thm_descent} ,\ref{thm_prod} and
 \cite[Proposition~3.2 (4)]{cao18approximation}.
}

\section*{Acknowledgment}
The authors would like to thank Ye Tian, Weizhe Zheng, Dasheng Wei,
 Yang Cao, Lei Zhang,
 Shizhang Li and Junchao Shentu for
 helpful discussions
  and the referees for valuable suggestions.
The work was partially done during the authors' visits to
 the Morningside Center of Mathematics, Chinese
 Academy of Sciences. They  thank the Center for its hospitality.

In preparing the second revision, we used OpenAI's ChatGPT/Codex
to assist with mathematical discussions, checking arguments and references,
revising the exposition, and drafting responses to the referees.
The authors retain full responsibility for the content.

\bibliography{unibib}
\bibliographystyle{amsalpha}
\end{document}